\documentclass[10pt,oneside]{article}
\usepackage{graphics}
\usepackage{amsmath}
\usepackage{amssymb}
\usepackage{enumerate}
\usepackage{mathrsfs}
\usepackage[dvips]{graphicx}
\usepackage{psfrag}
\usepackage{sidecap,epstopdf}
\usepackage{subfigure}
\usepackage{color}
\newtheorem{thm}{Theorem}
\newtheorem{lem}{Lemma}

\newcommand{\X}{\mathbf{x}}

\newcommand{\Dn}[2]{\frac{\partial #1}{\partial n_{#2}}}

\newcommand{\jump}[2]{\mu_O\Dn{#1}{}(\Bx)\Big|_{\partial \omega_\varepsilon^{(#2)+}}-\mu_{I_{#2}}\Dn{#1}{}(\Bx)\Big|_{\partial \omega_\varepsilon^{(#2)-}}}

\newcommand{\Oj}{\mathbf{O}}

\newcommand{\N}{\boldsymbol{n}}

\renewcommand{\theequation}{\thesection.\arabic{equation}}
\numberwithin{equation}{section}
\DeclareMathOperator*{\diag}{diag}

\title{Asymptotic analysis of solutions to transmission problems in solids with many inclusions}
\author{M.J. Nieves\footnote{Mechanical Engineering and Materials Research Centre, Liverpool John Moores University, James Parsons Building, Byrom Street, Liverpool L3 3AF, U.K.}}
\date{}

\newcommand{\bfm}[1]{\mbox{\boldmath ${#1}$}}

\newcommand{\beqa}{\begin{eqnarray}}
\newcommand{\eeqa}{\end{eqnarray}}
\newcommand{\bequ}{\begin{equation}}
\newcommand{\eequ}[1]{\label{#1}\end{equation}}

\newcommand{\GO}{\Omega}

\newcommand{\BGx}{\bfm\xi}

\newcommand{\BGL}{\bfm\Lambda}

\newcommand{\CP}{{\cal P}}

\newcommand{\CT}{{\cal T}}
\newcommand{\CU}{{\cal U}}
\newcommand{\CV}{{\cal V}}

\newcommand{\BCC}{{\bfm{\cal C}}}
\newcommand{\BCD}{{\bfm{\cal D}}}

\newcommand{\BCP}{{\bfm{\cal P}}}
\newcommand{\BCQ}{{\bfm{\cal Q}}}

\def\Bn{{\bf n}}

\def\Bw{{\bf w}}
\def\Bx{{\bf x}}
\def\By{{\bf y}}
\def\Bz{{\bf z}}

\def\BC{{\bf C}}

\def\BI{{\bf I}}

\def\BO{{\bf O}}
\def\BP{{\bf P}}
\def\BQ{{\bf Q}}

\def\BT{{\bf T}}

\def\BW{{\bf W}}
\def\BX{{\bf X}}
\def\BY{{\bf Y}}
\def\BZ{{\bf Z}}

\newcommand{\beq}{\begin{equation}}
\newcommand{\eeq}{\end{equation}}
\newcommand{\overliner}{\begin{eqnarray}}
\newcommand{\earr}{\end{eqnarray}}
\newcommand{\beqn}{\begin{equation*}}
\newcommand{\eeqn}{\end{equation*}}
\newcommand{\overlinern}{\begin{eqnarray*}}
\newcommand{\earrn}{\end{eqnarray*}}

\begin{document}
\maketitle
\begin{abstract}
We construct an asymptotic approximation to the solution of  
a  transmission problem for a body containing a region occupied by  many small inclusions. The cluster of inclusions is characterised by two small parameters that determine the nominal diameter of individual inclusions and their separation within the cluster. These small parameters can be comparable to each other. Remainder estimates of the asymptotic approximation  are rigorously justified. Numerical illustrations demonstrate the efficiency of the asymptotic approach when compared with benchmark finite element algorithms.
\end{abstract}
\section{Introduction}\label{intromesotran}

Uniform asymptotic approximations for solutions to boundary value problems involving large clusters of small defects have been constructed in the articles \cite{MM_meso_1, Maz_MMS, MMN_Mesoelast, MMN_Mesoelast_voids} and the monograph \cite{MMN_book}. The approximations have been developed for different operators of mathematical physics and for a range of different boundary conditions imposed on the surfaces of the small defects.  The approach employed in achieving these approximations does not utilise any strong assumptions on the arrangements of inclusions within the cluster, such as periodically distributed defects or arrangements which are statistically determined, where alternative conventional  techniques such homogenisation are applicable \cite{Bakhvalov, MarKhrus}.

Here we address the approximation of the solution to a transmission problem for a solid containing a dense non-periodic arrangement of small inclusions. In particular, such an approximation is capable of tracing the interaction of defects within a cluster, which is a serious challenge, especially in regions where fields are likely to  rapidly oscillate.

Several important approximations for solids containing dilute arrangements of defects, amongst much else,  have appeared in \cite{OPTH1, OPTH2}. There, the method of compound asymptotic expansions is  systematically presented for various elliptic boundary value problems of mathematical physics in singularly perturbed domains. For domains with small defects, the method relies on  model problems posed in domains without any holes, and problems in the exterior of individual small defects. The approach has also led to approximations for energy characteristics associated with these singularly perturbed problems in perforated domains such as eigenvalues, stress-intensity factors and capacities.

The method of compound asymptotic expansions has recently played a major role in  uniform asymptotic approximations for Green's kernels  in domains with several defects for both scalar  \cite{CRM, JCOM,  MMAS, Sob_vol} and vector problems \cite{RAN, AA}. In particular, a uniform approximations to singular fields for transmission problems in planar bodies containing several small inclusions  has appeared in \cite{AMS_tran}. Approximations of this type have been shown to provide results that give excellent comparison with those based on benchmark finite element schemes \cite{RAN, AMS_tran}. Uniform asymptotic approximations for Green's functions have also been used to model the flow of obstacles in Hele-Shaw flow \cite{Mishetal1,  Pecketal}.

The approximations mentioned above for dilute composites serve the case when the number of small defects are finite and are situated far apart from each other. However, in the situation when the number of defects becomes large and can be close to each other, one needs an alternative tool to model this scenario.

If the  small defects are arranged periodically in some region, then one can employ powerful homogenisation based techniques to model an effective medium \cite{MarKhrus}. The technique can also reveal additional contributions to the physics of such problems when the number of small particles within the region increase, while the overall volume occupied within a region remains constant \cite{Murat, MarKhrus}.  Periodic  composite materials for both electromagnetism and elasticity have been modelled using the homogenisation approach in \cite{SP} and has been extended to  treat problems where different boundary conditions are supplied on neighboring defects in \cite{Jager}.
This averaging procedure has led to the effective properties  of cubically arranged homogeneous spherical inclusions in an ambient matrix in \cite{Kristensson} and 
for  homogenisation to periodic elastic media with jumps in the transmission conditions on interfaces of small inclusions, see \cite{Orlik}.

Other approaches  used to establish effective behaviour of composites include a potential approach  used in \cite{Keller} to determine the effective conductivity for dense arrays of perfectly conducting spheres or perfectly and non-perfectly conducting cylinders. In addition, a functional equation approach has been used to study the  
effective conductivity of doubly periodic systems  of  inclusions distributed within a matrix having non-ideal contact conditions  in two dimensions \cite{Castro}.

In addition, the homogenisation technique can treat composites where defect positions do not exhibit periodicity, but may be specified by some statistical law. For example, see  \cite{Berlyand}  for a problem of this type that  considers the potential for medium containing randomly distributed circular inclusions under ideal  contact conditions.  In \cite{LPI, LPII, PonteWillis},   homogenisation based approximations have been used to obtain estimates of effective moduli characterising the composites for both elastic and hyperelastic materials containing randonly distributed fibres or defects.  For the case of a Neo-Hookean material containing periodically placed fibers that is subjected to different loading conditions, see \cite{Brunetal}.


Naturally,  one can find many examples of densely perforated materials for which the position of the perforations are not  governed by periodicity or a statistical law. Hence a homogenisation approach is not applicable when modelling these materials..
 The method of  meso-scale asymptotic approximations was developed in \cite{MM_meso_1} to approximate potentials for a bodies containing  large non-periodic clusters of small defects, with rigid boundaries. Meso-scale approximations for  solids containing a cluster of voids has appeared  in \cite{Maz_MMS}. More recently, the meso-scale approach has been used to approximate solutions for problems of the Lam\'e system for three-dimensional solids with clouds of defects with rigid boundaries in \cite{MMN_Mesoelast} and when the traction-free conditions are supplied on small voids \cite{MMN_Mesoelast_voids}.

Low-frequency vibration problems for solids with arrays of small inclusions have also been addressed using a modification of the methods of compound and meso-scale asymptotic approximations.  Asymptotics of the first eigenvalue and corresponding eigenfunction for domains with a cloud of rigid inclusions have appeared in \cite{EigMMN}.
Applications of the method of meso-scale approximations have also appeared in \cite{ChallaSini2, ChallaSini3}, where the scattering problems for many small obstacles in the infinite space were considered.

Here we adapt the approaches of \cite{MMN_book} and \cite{OPTH1, OPTH2}  to develop the approximation of the solution to a transmission problem inside a body with many small arbitrary inclusions.
Before stating the main result of this article, we supply the details of the problem we intend to tackle. Here, $\Omega$ will denote a bounded subset of $\mathbb{R}^3$,  which we assume contains a material with shear modulus $\mu_O$ and has smooth boundary $\partial\Omega$. Let $\omega$ be a region of $\Omega$, with diameter $1$. Contained in $\omega$ will be  many small inclusions $\omega_\varepsilon^{(k)}$, $1 \le k \le  N$. The $k^{th}$ inclusion has centre $\BO^{(k)}$, a smooth  interface $\partial \omega_\varepsilon^{(k)}$, a normalized diameter which is characterized by the small dimensionless parameter $\varepsilon$, and  is also occupied by a material with shear modulus $\mu_{I_k}$,   $1 \le k \le  N$.   Another small non-dimensional parameter $d$, defined by 
\[d=\frac{1}{2} \min_{\substack{j\ne k\\ 1 \le j, k \le N}} |\BO^{(j)}-\BO^{(k)}|\;, \]
is used to illustrate the ``closeness'' of one inclusion to the other within the cloud $\omega$.  
Additional geometric constraints  on $\omega$ are then given by
\begin{equation*}\label{assumptions}
 \text{dist}(\cup_{j=1}^N \omega^{(j)}_\varepsilon, \partial \omega)=2d\quad \text{ and } \quad \text{dist}(\omega, \partial \Omega)=1\;.
\end{equation*}
The parameters $N$ and $d$ satisfy the inequality
 \[N \le \text{const } d^{-3}\;.\]

Our main objective is to derive the asymptotic approximation to the displacement function $u_N$ satisfying the  transmission problem 
 
\begin{equation}\label{mesotranprob1a}
\left.\begin{array}{c}
\displaystyle{\mu_O \Delta u_N(\Bx)=f(\Bx)\;, \quad \Bx \in \Omega_N=\Omega \backslash \cup_{ k=1}^N \overline{\omega_\varepsilon^{(k)}}\;,}
\\\\
\displaystyle{\mu_{I_j} \Delta u_N(\Bx)=0\;, \quad \Bx \in \omega_\varepsilon^{(j)}, j=1, \dots, N\;,}
\\\\
\displaystyle{u_N(\Bx)=\phi(\Bx)\;, \quad \Bx\in \partial \Omega\;,}
\\ \\
\displaystyle{u_N(\Bx)\Big|_{\partial \omega^{(j)+}_\varepsilon}=u_N(\Bx)\Big|_{\partial \omega^{(j)-}_\varepsilon}\;, \quad j=1, \dots, N\;,}
\\ \\
\displaystyle{\mu_O\Dn{u_N}{}(\Bx)\Big|_{\partial \omega_\varepsilon^{(j)+}}=\mu_{I_j}\Dn{u_N}{}(\Bx)\Big|_{\partial \omega_\varepsilon^{(j)-}}\;, \quad j=1, \dots, N\;,}
\end{array}\right\}
\end{equation}
where $\partial \omega_\varepsilon^{(j)+}$ ($\partial \omega_\varepsilon^{(j)-}$) represents the boundary $\partial \omega^{(j)}_\varepsilon$ approached from the exterior (interior).
The function $f$,describing the body force in $\Omega_N$,  belongs to the space $L_2(\Omega)$,  and  has a support satisfying $ \omega\cap\text{supp }f  =\varnothing$ and $\text{dist}(\omega, \text{supp }f)=O(1)$. In the displacement condition on $\partial \Omega$ we have $\phi \in L^{1/2, 2}(\partial \Omega)$.

The   construction of the asymptotics of   $u_N$ relies on the methods of compound and meso-scale  asymptotic expansions, which in turn makes use of model fields defined in the unperturbed  set (without small inclusions) $\Omega$  and in the infinite space containing a single small inclusion $\omega_\varepsilon^{(k)}$, $k=1, \dots, N$. Such model fields involve:
\begin{enumerate}
\item the solution  $w_f$ of the Dirichlet problem  of Poisson's equation in $\Omega$; 
\item the vector functions $\BCD^{(k)}_\varepsilon$,     whose components are the dipole fields for the inclusion $\omega_\varepsilon^{(k)}$. These fields allow one to construct  boundary layers outside small holes in the asymptotic algorithm;
\item the regular part  $H$ of Green's function $G$ in $\Omega$.
\end{enumerate}
It will also be shown that coefficients near boundary layers in the approximation of $u_N$ form solutions to a  certain algebraic system. This system involves derivatives of $w_f$ and  integral characteristics pertaining to the small voids 
such as the polarization matrix $\BCP^{(k)}_\varepsilon$, (see 
 \cite{MMP})  which  is a    $3\times 3$  matrix for the small inclusions $\omega_\varepsilon^{(k)}$, $1\le k \le N$.  As is discussed in more detail  later, this matrix can be positive or negative definite.
In addition, if $\BCP^{(k)}_\varepsilon$ is negative (positive) definite, we assume that the  maximum and minimum eigenvalues $\lambda^{(j)}_{\text{max}}$ and $\lambda^{(j)}_{\text{min}}$, respectively,  of $-\BCP^{(j)}_\varepsilon$  ($\BCP^{(j)}_\varepsilon)$ satisfy
\begin{equation}\label{eigest}
C_1\varepsilon^3 <\lambda_{min}^{(j)},\qquad \text{ and }\qquad  \lambda_{max}^{(j)}<  C_2\,\varepsilon^3,
\end{equation}
where $C_1$ and $C_2$ are constants independent of $\varepsilon$.


\vspace{0.1in}
\begin{thm}
 \label{thm1_alg_f_trans} 
Let  
\begin{equation}\label{epscd}
\varepsilon < c\, d\;,
\end{equation}
where $c$ is a sufficiently small absolute constant. Then the solution $u_N(\Bx)$ 
admits the asymptotic representation
\begin{equation}\label{introeq1}
u_N(\Bx)=w_f
(\Bx)+
 \sum_{1\le k\le N} \BC^{(k)}\cdot \{ \BCD^{(k)}_\varepsilon(\Bx)-\BCP^{(k)}_\varepsilon \nabla_\By H(\Bx, \By)\Big|_{\By=\Oj^{(k)}}\} 
+R_N(\Bx)\;,
\end{equation}
where $\BC^{(k)}=(C^{(k)}_1, C^{(k)}_2, C^{(k)}_3)^T$, $1\le k \le N$ 
satisfy the solvable linear
algebraic system 
\beq
 \nabla w_f(\BO^{(j)})+\BC^{(j)}
+\sum_{\substack{k \ne j\\ 1\le k\le N}} (\nabla_{\Bz}\otimes \nabla_{\Bw})G(\Bz, \Bw)\Big|_{\substack{\Bz=\Oj^{(j)}\\ \Bw=\BO^{(k)}}}\BCP_\varepsilon^{(k)}\BC^{(k)}=\BO\;,\quad \text{ for  }j=1, \dots, N\;.
\eequ{alg_s_intro}
The 
remainder $R_N$ satisfies the energy
estimate
\begin{equation}\label{introeq2}
  \int_{\cup_{k=1}^N \omega^{(k)}_\varepsilon\cup \Omega_N} |\nabla R_N|^2\, d\Bx
 \le \text{\emph{const} } \Big\{ \varepsilon^{11
 }d^{-11
 } + \varepsilon^{5}d^{-3} \Big\} \| \nabla w_f
 \|^2_{L_2( \GO 
 )}  
 \end{equation}
\end{thm}

\begin{figure}\centering
     \subfigure[][]{
\centering
        \includegraphics[width=0.45\textwidth]{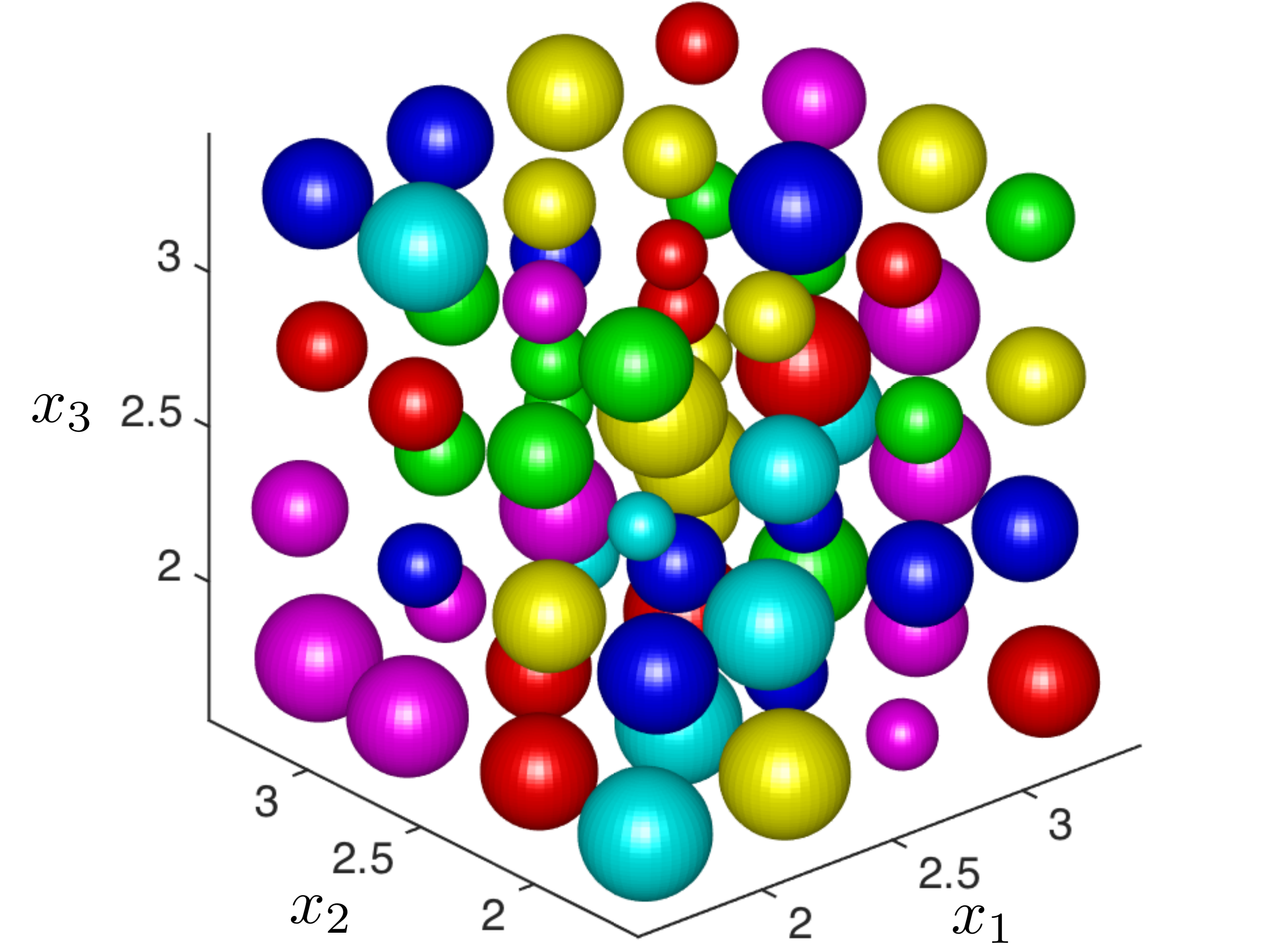}
          \label{fig:3ab}}
          \\
         \subfigure[][]{
\centering
        \includegraphics[width=0.45\textwidth]{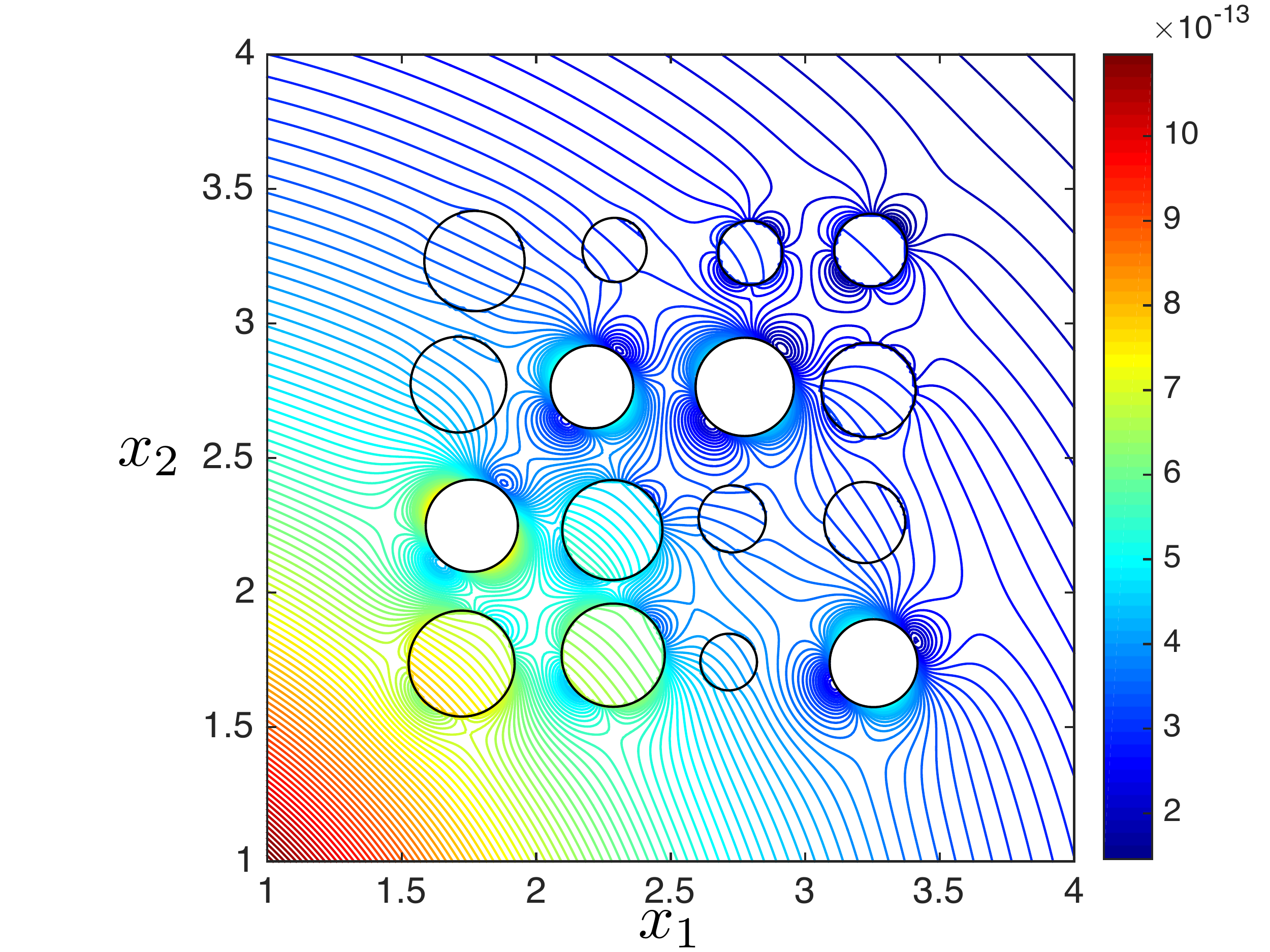}\label{fig:4a}}
           \subfigure[][]{
\centering
        \includegraphics[width=0.45\textwidth]{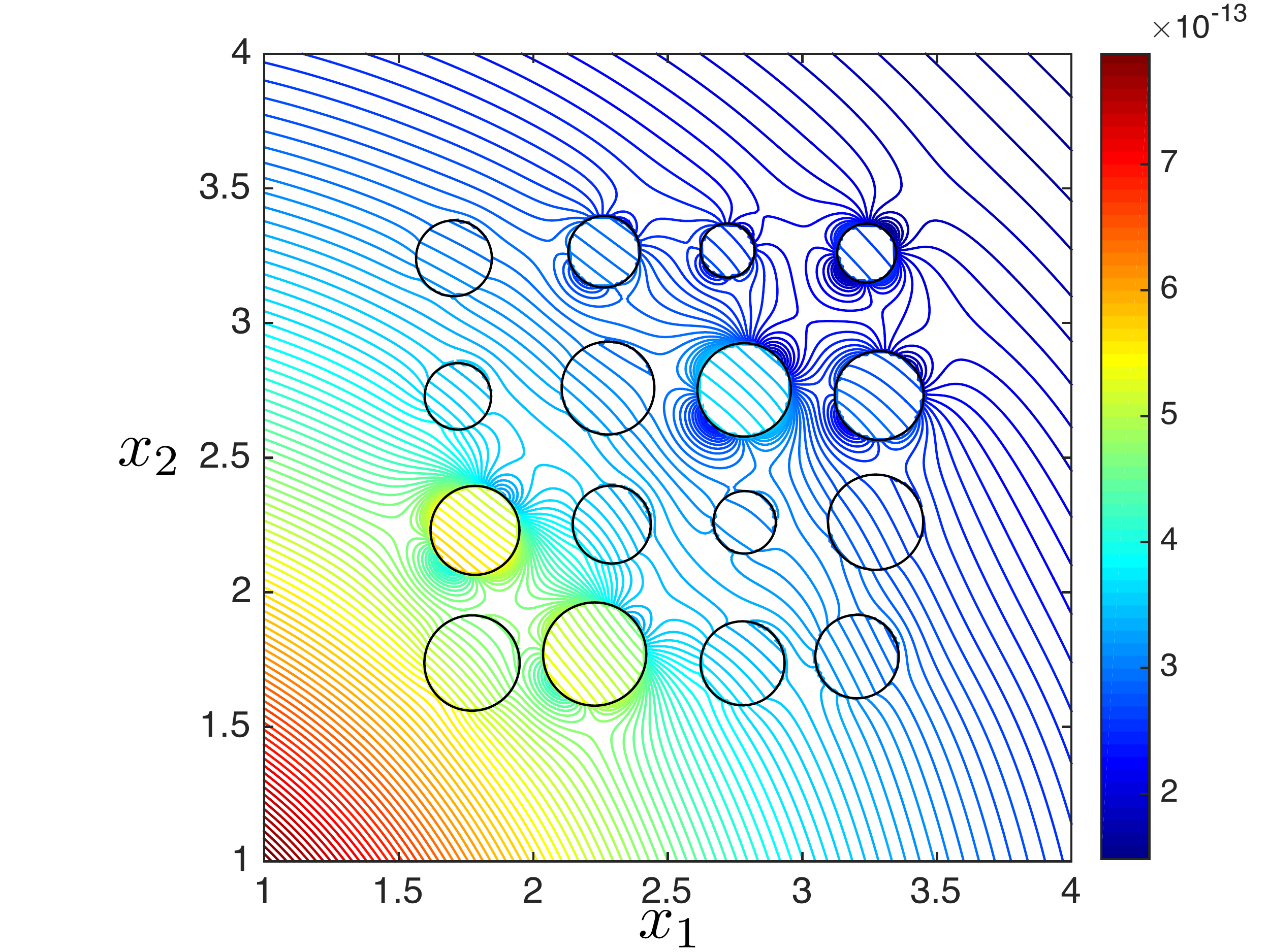}\label{fig:4b}}\\
        \subfigure[][]{
\centering
        \includegraphics[width=0.45\textwidth]{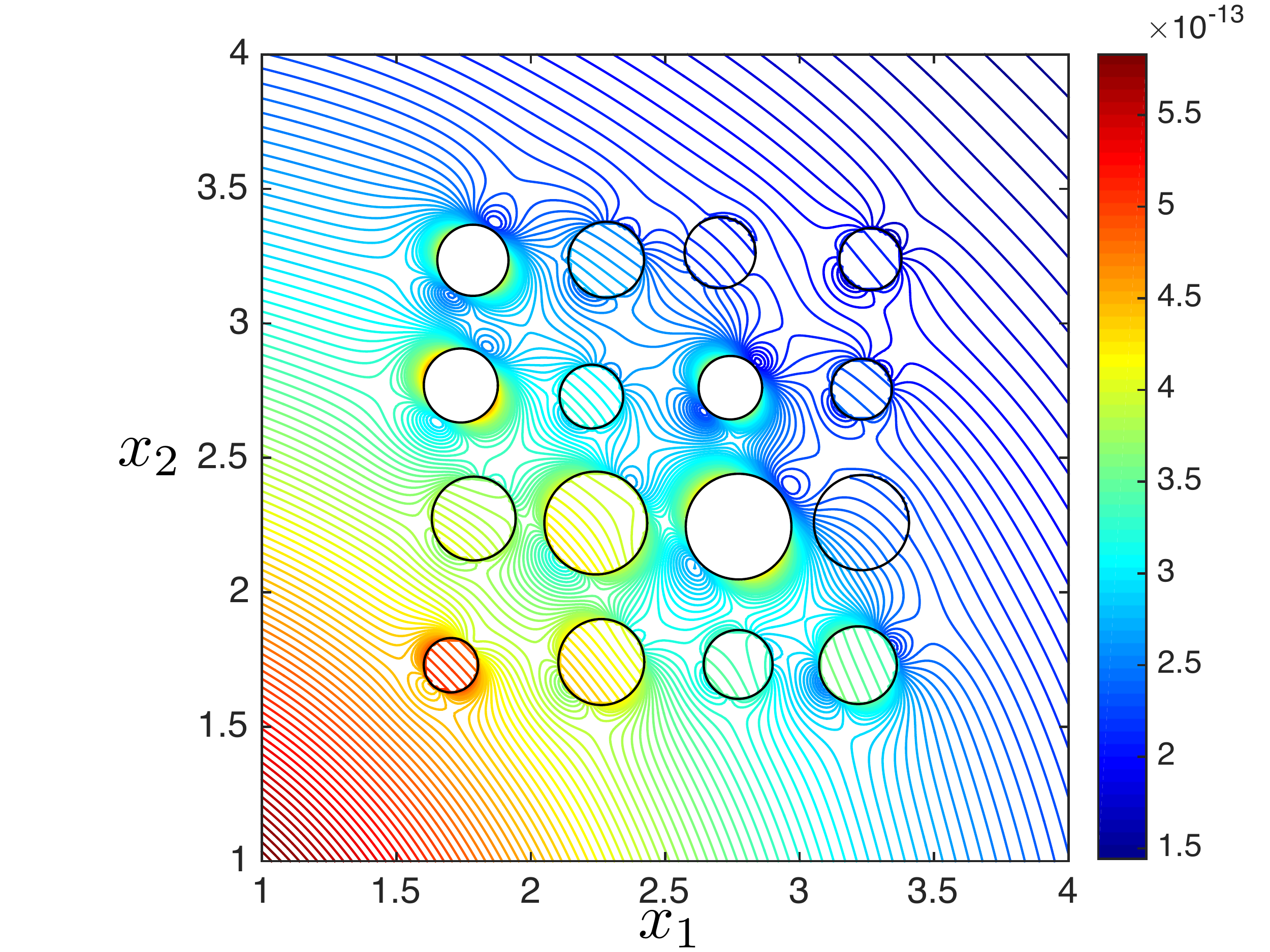}\label{fig:5a}}
           \subfigure[][]{
\centering
        \includegraphics[width=0.45\textwidth]{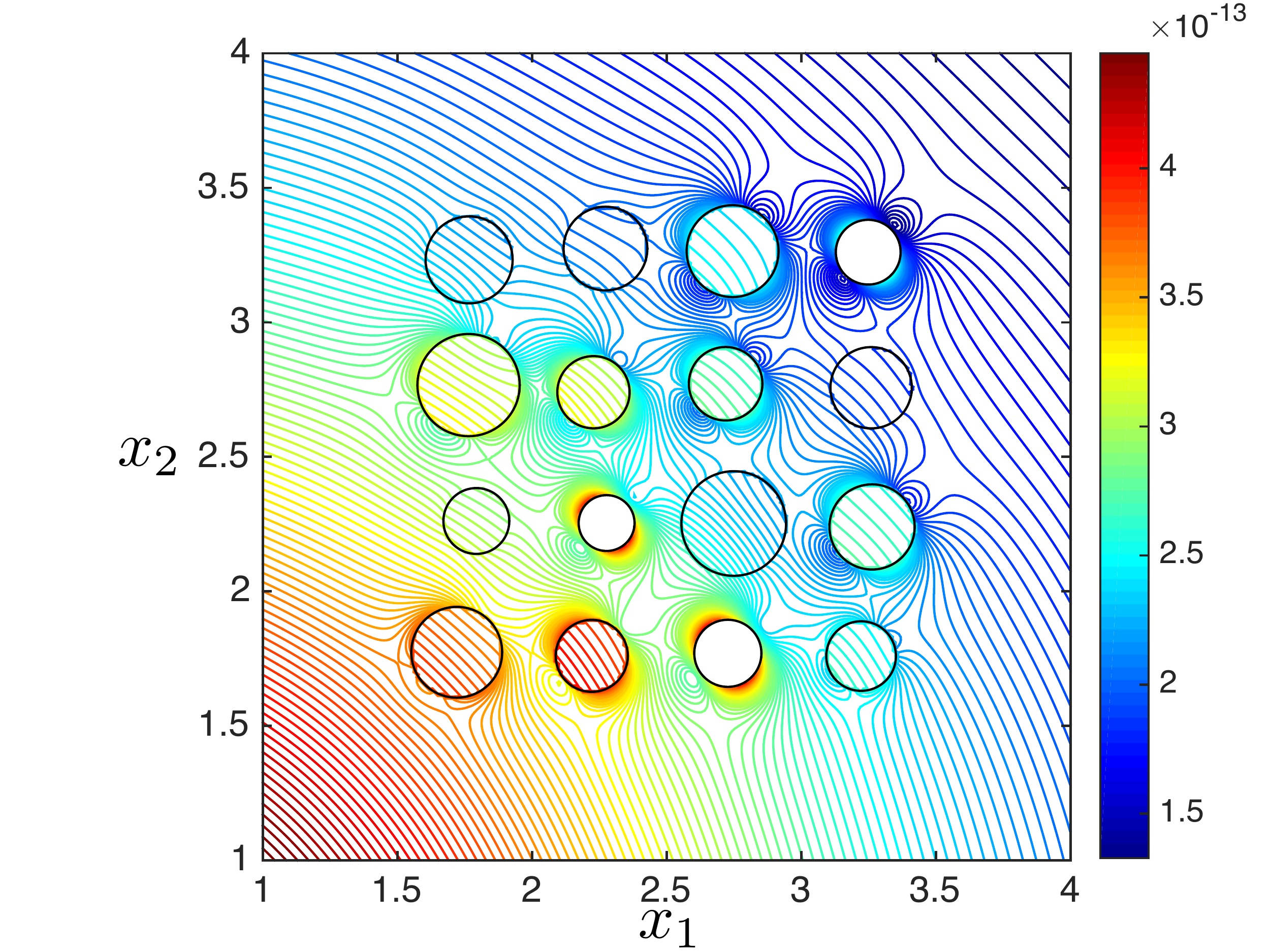}\label{fig:5ba}}

\caption[]{(a) A cluster of 64 inclusions. The colors shown indicate the material contained in the inclusion: Cast Iron (green), Steel AISI 4340 (blue), Aluminum (yellow), Copper (light blue) and Iron (purple). We assume the ambient matrix is occupied by Structural Steel. Here, the inclusions which are red  correspond to voids (which are not occupied by a material). Computations for $|\nabla u_N|$, based on the asymptotic approximation  (\ref{introeq1}), in the vicinity of the cluster are shown in (b)--(d)  along the cut-planes (b) $x_3=1.75$, (c) $x_3=2.25$, (d) $x_3=2.75$ and (d)  $x_3=3.25$.}
\label{fig:3}%
\end{figure}

As an example, we consider a large cluster of  inclusions ($N=64$) arranged in cube, as according to Figure \ref{fig:3}, which is embedded in sphere of radius 7.  Here $f(\Bx)$ in (\ref{mesotranprob1a}) is a radially symmetric function having support inside the ball of radius 1.5 (further details of the numerical set up can be found in section \ref{numericalsimulations}). The cluster is composed of both voids and inclusions which are occupied by either Cast Iron, Steel AISI 4340, Aluminum,  Copper or Iron. The ambient matrix has the material properties of Structural Steel. For such a problem, the task of using the method of finite elements, with the package COMSOL, to obtain the solution $u_N$   can be computationally intensive. In fact, COMSOL could not compute the solution to this problem in this case. However, the asymptotic formulae (\ref{introeq1}) remains efficient and we present computations for $|\nabla u_N|$, based on the derivatives of the leading order approximation to $u_N$ in (\ref{introeq1}), along cut-planes which intersect the cloud.

The structure of the article is as follows. In section \ref{modprob_incl} we introduce model problems which allow one to construct the formal approximation to $u_N$, and this formal asymptotic procedure is provided in section \ref{Formalgtran}. There, the algebraic system, concerning coefficients involved in the asymptotic approximation to $u_N$, is determined and the solvability of this system is discussed in section \ref{algmesotran}. The proof of the energy estimate (\ref{introeq2}) for the remainder term involved in the approximation  is provided in section \ref{energymesotran}, where the completion of the proof of Theorem \ref{thm1_alg_f_trans} is also presented. In section \ref{corosec}, we extend the results of Theorem \ref{thm1_alg_f_trans} to the case of the transmission problem for an infinite medium containing a non-periodic cluster of small inclusions. Following this, in section \ref{connection}, we further investigate the algebraic system (\ref{alg_s_intro}) for a large periodic cluster situated inside a body and from this we derive an auxiliary homogenised problem concerning the effective inclusion $\omega$ situated inside the domain $\Omega$. Numerical illustrations are then given in section \ref{numericalsimulations} that show the efficiency of the asymptotic approach, in particular, when compared with computations based on finite element algorithms.
In section \ref{conclusions}, we give some conclusions and discussion. Finally, in the Appendix, we present the proofs of auxiliary results needed to prove the solvability of the algebriac system discussed in section \ref{algmesotran}.
\section{Model problems} \label{modprob_incl}
We now introduce model problems posed in either the unperturbed domain $\Omega$ or the infinite space with the small inclusion at the origin, which we use in the formal asymptotic procedure given in the following section.

\begin{enumerate}[$1.$]
\item \emph{The solution $w_f$ of the Dirichlet problem for Poisson's equation. }Let $w_f$ satisfy the problem
\begin{equation}\label{goveqw_f}
\mu_O\Delta w_f (\Bx)=f(\Bx)\;, \quad \Bx \in \Omega\;,
\end{equation}
\[w_f(\Bx)=\phi(\Bx)\;, \quad \Bx \in \partial \Omega\;,\]
where as before $f \in L_2(\Omega)$, $\text{supp }f \cap \omega=\varnothing$, $\text{dist}(\omega, \text{supp }f)=O(1)$ and $\phi \in L^{1/2, 2}(\partial \Omega)$. Later, we assume that $f$ is extended by zero inside the inclusions $\omega^{(k)}_\varepsilon$, $1\le k \le N$.

\item \emph{The regular part $H$ in $\Omega$. }Let $H$ be the regular part of Green's function $G$ in $\Omega$, which satisfies
\[\mu_O\Delta_\Bx H(\Bx, \By)=0\;, \quad \Bx, \By \in \Omega\;,\]
\begin{equation}\label{bcH}
H(\Bx, \By)=(4\pi\mu_O|\Bx-\By|)^{-1}\;, \quad \Bx \in \partial \Omega\;, \By \in \Omega\;,
\end{equation}
and $G$ is related to $H$ via
\begin{equation}\label{Grep}
H(\Bx, \By)=\frac{1}{4\pi \mu_O|\Bx-\By|}-G(\Bx, \By)\;.
\end{equation}
  \item \emph{The dipole fields for the small inclusion $\omega_\varepsilon^{(j)}$, $j=1, \dots, N$.} The dipole fields $\mathcal{D}^{(k)}_j$, $j=1, 2, 3$, for the scaled inclusion $\omega^{(k)}=\{\BGx: \varepsilon \BGx+\BO^{(j)}\in \omega_\varepsilon^{(k)}\}$, are now introduced as components of the vector function $\BCD^{(k)}=(\mathcal{D}^{(k)}_1, \mathcal{D}^{(k)}_2, \mathcal{D}^{(k)}_3)^T$, $k=1, \dots, N$, which solves the transmission problem
  \begin{equation}\label{dipprob}
  \left.\begin{array}{c}
  \mu_O\displaystyle{\Delta \BCD^{(k)}(\BGx)=\BO\;, \quad \BGx \in \mathbb{R}^3 \backslash \overline{\omega^{(k)}}\;,}\\
  \\
 \mu_{I_k}\displaystyle{ \Delta \BCD^{(k)}(\BGx)=\BO\;, \quad \BGx \in \omega^{(k)}\;,}\\
  \\
  \displaystyle{\BCD^{(k)}(\BGx)\Big|_{\partial \omega^{(k)+}}=\BCD^{(k)}(\BGx)\Big|_{\partial \omega^{(k)-}}\;, }
\\ \\
\displaystyle{\mu_O\Dn{\BCD^{(k)}}{}(\BGx)\Big|_{\partial \omega^{(k)+}}-\mu_{I_k}\Dn{\BCD^{(k)}}{}(\BGx)\Big|_{\partial \omega^{(k)-}}=(\mu_O-\mu_{I_k})\N^{(k)}\;, }
   \end{array}\right\}
  \end{equation}
 where $\N^{(k)}$ is the unit outward normal to $\omega^{(k)}$.
  The leading order behaviour of $\BCD^{(k)}$, described in (\ref{dipprob}),  can be written explicitly using the $3\times 3$ symmetric polarization matrix $\BCP^{(k)}=\{\CP_{ij}^{(k)}\}_{i, j=1}^3$ for the small inclusion. We have
  
  \begin{lem}\label{lemdipasymp} $($\emph{see \cite{MMP}}$)$
Let $|\BGx|>2$, then the vector function $\BCD^{(k)}$ admits the asymptotic representation
  \begin{equation}\label{asydipeq}
  \BCD^{(k)}(\BGx)= -\BCP^{(k)}\nabla ((4\pi\mu_O|\BGx|)^{-1})+O(|\BGx|^{-3})\;,
  \end{equation}
  where the entries of $\BCP^{(k)}=[\CP^{(k)}_{ij}]_{i,j=1}^3$ are given by 
  \begin{eqnarray}
  \CP^{(k)}_{ij}&=&(\mu_{I_k}-\mu_O)\text{{\rm meas\rm}}(\omega^{(k)})\delta_{ij}-\mu_{O}\int_{\mathbb{R}^3\backslash \overline{\omega^{(k)}}} \nabla\mathcal{D}^{(k)}(\BGx)\cdot \nabla\mathcal{D}^{(j)}(\BGx) d\BGx\nonumber \\
  &&-\mu_{I_k}\int_{\omega^{(k)}} \nabla\mathcal{D}^{(k)}(\BGx)\cdot \nabla\mathcal{D}^{(j)}(\BGx) d\BGx\;,\label{expressionP}
  \end{eqnarray}
  for $1\le i,j \le 3$.
 \end{lem}
\end{enumerate}

Here, the expression (\ref{expressionP}) shows that the  polarization matrix $\BCP^{(k)}$ is symmetric. Further, after integration by parts and using the definition (\ref{dipprob}) of the components of the vector function $\BCD$ one can show that
\begin{equation}\label{CPG}
\CP^{(k)}_{ij}=(\mu_{I_k}-\mu_O)\Big\{\text{{\rm meas\rm}}(\omega^{(k)})\delta_{ij}+\int_{\partial \omega^{(k)} }\mathcal{D}_i(\BGx) \frac{\partial \BGx_j}{\partial n}dS_{{ \BGx}}\Big\}\;,
\end{equation}
where further integration by parts in $\omega^{(k)}$ in the term on the right gives
\[\int_{\partial \omega^{(k)} }\mathcal{D}_i(\BGx)\Big|_{\partial \omega^{(k)-}} \frac{\partial \BGx_j}{\partial n}dS_{\BGx}=\int_{\partial \omega^{(k)} }\BGx_j \frac{\partial  \mathcal{D}_i(\BGx)}{\partial n}\Big|_{\partial \omega^{(k)-}}dS_{{\BGx}}\;.\]
As a result the integral term in (\ref{CPG}) defines a Gram matrix.
 Here (\ref{CPG}) shows that for voids ($\mu_{I_k}=0$)  the dipole matrix is negative definite, whereas if $\mu_{I_k}>\mu_O$ ($\mu_{I_k}<\mu_O$), this polarization tensor is  positive (negative) definite. 

By then rescaling, we  create the vector functions $\BCD^{(j)}_\varepsilon(\Bx)=\varepsilon\BCD(\varepsilon^{-1}(\Bx-\BO^{(j)}))$ and matrices $\BCP^{(j)}_\varepsilon=\varepsilon^3\BCP^{(j)}$, $1\le j\le N$, which are to be used throughout the asymptotic algorithm. Here, the components of  $\BCD^{(j)}_\varepsilon$ are then the dipole fields for the inclusion $\omega^{(k)}_\varepsilon$, $1\le k \le N$.


\section{Formal asymptotic procedure}\label{Formalgtran}
We now formally construct an asymptotic approximation to the solution $u_N$ of the transmission problem (\ref{mesotranprob1a}).
We prove the following Lemma.
\begin{lem}\label{lemformapp}
The formal asymptotic approximation of $u_N$, the solution of $(\ref{mesotranprob1a})$, has the form
\begin{equation*}\label{formalapp}
u_N(\Bx)=w_f(\Bx)+\sum_{1 \le k \le N} \BC^{(k)} \cdot \{ \BCD^{(k)}_\varepsilon(\Bx)-\BCP_\varepsilon^{(k)} \nabla_{\By} H(\Bx, \By)\Big|_{\By=\BO^{(k)}}\}+R_N(\Bx)\;,
\end{equation*}
where $\BC^{(k)}$, $1\le k \le N$, are solutions of the algebraic system
\begin{equation*}\label{indivualalgeq_1}
 \nabla w_f(\BO^{(j)})+\BC^{(j)}
+\sum_{\substack{k \ne j\\ 1\le k\le N}} (\nabla_{\Bz}\otimes \nabla_{\Bw})G(\Bz, \Bw)\Big|_{\substack{\Bz=\Oj^{(j)}\\ \Bw=\BO^{(k)}}}\BCP_\varepsilon^{(k)}\BC^{(k)}=\BO\;,\quad \text{ for  }j=1, \dots, N\;.
\end{equation*}
  The remainder term $R_N$ then satisfies the problem
 \begin{eqnarray*}
 &&\mu_O \Delta R_N(\Bx)=0\;, \quad \Bx\in \Omega\backslash \cup_{j=1}^N \overline{\omega^{(j)}_\varepsilon}, \quad \mu_{I_k} \Delta R_N(\Bx)=0\;, \quad \Bx \in \cup_{j=1}^N {\omega^{(j)}_\varepsilon},\\
 && R_N(\Bx)=\psi(\Bx)\;, \quad \Bx\in \partial \Omega\\
 &&R_N(\Bx)\Big|_{\partial \omega^{(k)+}_\varepsilon}=R_N(\Bx)\Big|_{\partial \omega^{(k)-}_\varepsilon}\\
 &&\frac{\partial R_N}{\partial n}(\Bx)\Big|_{\partial \omega^{(k)+}_\varepsilon}-\frac{\partial R_N}{\partial n}(\Bx)\Big|_{\partial \omega^{(k)-}_\varepsilon}=\varphi^{(j)}_\varepsilon(\Bx)\;.
  \end{eqnarray*} 
 where 
 \begin{eqnarray*}
 &&|\psi(\Bx)|=O\Big(\sum_{1\le k\le N} \varepsilon^4 |\BC^{(k)}| |\Bx-\BO^{(k)}|^{-3}\Big)\\
 &&|\varphi^{(j)}_\varepsilon(\Bx)|=O\Big(\varepsilon\Big\{1+\varepsilon^2 |\BC^{(j)}|+\sum_{\substack{k \ne j\\ 1 \le k \le N}} \frac{\varepsilon^3 |\BC^{(k)}|}{|\BO^{(j)}-\BO^{(k)}|^4}\Big\}\Big)\;.
  \end{eqnarray*} 
 \end{lem}

{\emph{Proof. }}We first attempt to satisfy the governing equations in $\bigcup_{j=1}^N \omega_\varepsilon^{(j)} \cup \Omega_N$ and the exterior boundary condition appearing in (\ref{mesotranprob1a}).
Therefore, we approximate $u_N$ by the field $w_f$ defined in $\Omega$, i.e.
\begin{equation}\label{uNapp1}
u_N(\Bx)=w_f(\Bx)+W_N(\Bx)\;.
\end{equation}
Considering the boundary value problem for $W_N$, we have 
\[\mu_O \Delta W_N(\Bx)=0\;, \quad \Bx \in \Omega_N\;, \]
\[\mu_{I_j} \Delta W_N(\Bx)=0\;,\quad \Bx \in \omega_\varepsilon^{(j)}, j=1, \dots, N\;,\]
\[W_N(\Bx)=0\;, \quad \Bx \in \partial \Omega\;.\]
Since $w_f$ is defined everywhere in $\Omega$, $W_N$ is continuous across the frontier $\partial \omega_\varepsilon^{(j)}$, $j=1, \dots, N$. However,  the jump in the traction condition across $\partial \omega_\varepsilon^{(j)}$, $j=1, \dots, N$, is
\[\jump{{W_N}}{j}=-(\mu_O-\mu_{I_j}) \N^{(j)} \cdot \nabla w_f(\Bx)\;,\quad j=1, \dots, N\;.\]
Since the inclusion $\omega_\varepsilon^{(j)}$, $j=1, \dots, N$, is small, we can use  Taylor's expansion of  the derivatives of $w_f$ about $\Oj^{(j)}$ in the preceding condition, to obtain
\begin{equation}\label{tranalg1}
\jump{{W_N}}{j}=-(\mu_O-\mu_{I_j}) \N^{(j)} \cdot \nabla w_f(\Oj^{(j)})+O(\varepsilon)\;,\quad j=1, \dots, N\;.
\end{equation}
We now use the notion of the dipole fields for the inclusions (see Problem 3, Section \ref{modprob_incl}) to compensate for the error in the above right-hand side and construct $W_N$ with the representaion
\begin{equation}\label{uNapp2}
W_N(\Bx)=\sum_{1 \le k \le N} \BC^{(k)} \cdot \{ \BCD^{(k)}_\varepsilon(\Bx)-\BCP_\varepsilon^{(k)} \nabla_{\By} H(\Bx, \By)\Big|_{\By=\BO^{(k)}}\}+R_N(\Bx)\;,
\end{equation}
where in subsequent steps we will identify the algebraic system satidfied by $\BC^{(k)}=(C_1^{(k)}, C_2^{(k)}, C_3^{(k)})^T$, $k=1, \dots, N$.

Then, the field $R_N$ is harmonic inside $\Omega_N$ and $\omega_\varepsilon^{(j)}$, $j=1, \dots, N$. The asymptotics of $\BCD^{(k)}_\varepsilon$,  (see Lemma \ref{lemdipasymp}), and the boundary condition (\ref{bcH}) for $H$, allow one to assert that
\begin{equation}\label{uNapp3}
R_N(\Bx)=O\Big(\sum_{1\le k\le N} \varepsilon^4 |\BC^{(k)}| |\Bx-\BO^{(k)}|^{-3}\Big)\;, \quad \Bx \in \partial \Omega\;. \end{equation}
On the other hand, the same far-field representation for $\BCD_\varepsilon^{(k)}$, $k \ne j$, in (\ref{tranalg1}) provides the  displacement condition
\[R_N(\Bx)\Big|_{\partial \omega^{(j)+}_\varepsilon}=R_N(\Bx)\Big|_{\partial \omega^{(j)-}_\varepsilon}\;,\quad  1\le j \le N\;,\]
and the traction condition (see (\ref{tranalg1}))
\begin{eqnarray*}
&&\jump{R_N}{j}\\
&=&-(\mu_O-\mu_{I_j})\N^{(j)}\cdot \Big\{ \nabla w_f(\BO^{(j)})+\BC^{(j)}
+\sum_{\substack{k \ne j\\ 1\le k\le N}} (\nabla_{\Bz}\otimes \nabla_{\Bw})G(\Bz, \Bw)\Big|_{\substack{\!\!\!\Bz=\Bx\\ \Bw=\BO^{(k)}}}\BCP_\varepsilon^{(k)}\BC^{(k)}\\ 
&&+O(\varepsilon)+O(\varepsilon^3 |\BC^{(j)}|) +O\Big(\sum_{\substack{k \ne j\\ 1 \le k \le N}} \frac{\varepsilon^4 |\BC^{(k)}|}{|\Bx-\BO^{(k)}|^4}\Big)\Big\}, \quad j=1, \dots, N\;.
\end{eqnarray*}
Next, we expand the second order derivatives of $G$ about $\Bx=\BO^{(j)}$, to give
\begin{eqnarray}
&&\jump{R_N}{j}\nonumber \\
&=&-(\mu_O-\mu_{I_j})\N^{(j)}\cdot \Big\{ \nabla w_f(\BO^{(j)})+\BC^{(j)}
+\sum_{\substack{k \ne j\\ 1\le k\le N}} (\nabla_{\Bz}\otimes \nabla_{\Bw})G(\Bz, \Bw)\Big|_{\substack{\!\!\!\Bz=\BO^{(j)}\\ \Bw=\BO^{(k)}}}\BCP_\varepsilon^{(k)}\BC^{(k)}\nonumber \\ 
&&+O(\varepsilon)+O(\varepsilon^3 |\BC^{(j)}|) +O\Big(\sum_{\substack{k \ne j\\ 1 \le k \le N}} \frac{\varepsilon^4 |\BC^{(k)}|}{|\BO^{(j)}-\BO^{(k)}|^4}\Big)\Big\}, \quad j=1, \dots, N\;.\label{uNapp4}
\end{eqnarray}
Inspecting the last condition then suggests that $\BC^{(k)}$, $1\le j \le N$ should satisfy 
 \begin{equation}\label{indivualalgeq}
 \nabla w_f(\BO^{(j)})+\BC^{(j)}
+\sum_{\substack{k \ne j\\ 1\le k\le N}} (\nabla_{\Bz}\otimes \nabla_{\Bw})G(\Bz, \Bw)\Big|_{\substack{\Bz=\Oj^{(j)}\\ \Bw=\BO^{(k)}}}\BCP_\varepsilon^{(k)}\BC^{(k)}=\BO\;,\quad \text{ for  }j=1, \dots, N\;,
\end{equation}
to allow for the removal of the leading order discrepancy in the preceding traction condition.
Finally, combining (\ref{uNapp1}) and (\ref{uNapp2})--(\ref{indivualalgeq}) completes the proof of Lemma \ref{lemformapp}.
\hfill $\Box$
\section{Algebraic system and its solvability}\label{algmesotran}
Here we prove a result concerning the solvability of the algebraic system (\ref{indivualalgeq}).
\begin{lem}\label{lemsolvability}
Let 
\begin{equation}\label{eqepstod}
\varepsilon< \, c\, d\;,
\end{equation}
where $c$ is a sufficiently small absolute constant. Then the linear algebraic system $(\ref{indivualalgeq})$ is solvable and the estimate 
\begin{equation}\label{sumest1}
\sum_{j=1}^N |\BC^{(j)}|^2\le \text{\emph{Const} }\sum_{j=1}^N |\nabla u(\BO^{(j)})|^2\;,
\end{equation}
holds.
\end{lem}
We postpone the proof of this Lemma, in order to rewrite the algebraic system in a way which will simplify its representation, and we give  some auxiliary results.
\subsubsection*{Representation of the algebraic system and auxiliary results}
We begin by rewriting this system as follows. Set 
\[\BGL=((\nabla w_f(\BO^{(1)}))^T, \dots, (\nabla w_f(\BO^{(N)}))^T)\]
and 
\[\BCC=((\BC^{(1)})^T, \dots, (\BC^{(N)})^T)^T\;.\]
Next, we define $\BP_\varepsilon$ to be a $3N \times 3N$ block diagonal matrix  given by
\[\BP_\varepsilon=\diag \{ \BCP^{(1)}_\varepsilon, \dots, \BCP^{(N)}_\varepsilon\}\;,\]
and let  $\BT=[T_{ij}]_{i, j=1}^N$,  be another $3N \times 3N$ matrix, with $3 \times 3$ block entries $T_{ij}$ represented by
\[T_{ij} =\left\{\begin{array}{ll}
(\nabla_{\Bz} \otimes \nabla_{\Bw}) G(\Bz, \Bw)\Big|_{\substack{\Bz=\BO^{(j)}\\ \Bw=\BO^{(k)}}}\;, &\quad \text{ when }j \ne k\;,\\
 0I_3\;, &\quad \text{otherwise}\;.
\end{array}\right.\]
Then, (\ref{indivualalgeq}) takes the equivalent form:
\begin{equation}\label{system2}
\BCC+\BT \BP_\varepsilon \BCC=-\BGL\;.
\end{equation}
In addition, we introduce the matrix $\BQ_\varepsilon=\text{diag}\{\BCQ^{(1)}_\varepsilon,\dots, \BCQ^{(N)}_\varepsilon\}$ which is a $3N\times 3N$ block diagonal matrix where
\begin{equation}
\BCQ_\varepsilon^{(j)} =\left\{\begin{array}{ll}
-\BCP_\varepsilon^{(j)}\;, &\quad \text{ if  $\BCP^{(j)}_\varepsilon$ is negative definite}\;,\\
 \BCP_\varepsilon^{(j)}\;, &\quad \text{ if  $\BCP^{(j)}_\varepsilon$ is positive definite}\;.
\end{array}\right.\label{eqDeps}
\end{equation}

Finally, before  presenting the proof of Lemma \ref{lemsolvability}, we note the following result.
\begin{lem}\label{lemest_innerprod}
The estimate 
\[|\langle\BT\BP_\varepsilon \BCC,\BQ_\varepsilon \BCC\rangle| \le \text{\emph{Const} } d^{-3}\langle\BQ_\varepsilon \BCC, \BQ_\varepsilon \BCC\rangle \;.\]
holds.
\end{lem}
The  proof of the preceding Lemma is found in the Appendix.

\subsubsection*{Proof of Lemma \ref{lemsolvability}}
Taking the scalar product of (\ref{system2}) with $\BQ_\varepsilon \BCC$ we arrive at 
\begin{equation}\label{scalareq}
\langle\BCC, \BQ_\varepsilon \BCC\rangle+\langle\BT\BP_\varepsilon \BCC,\BQ_\varepsilon \BCC\rangle=-\langle\BGL,\BQ_\varepsilon \BCC \rangle\;.
\end{equation}
We apply the Cauchy inequality to the right-hand side to get 
\begin{equation*}\label{scalareq1}
\langle\BCC, \BQ_\varepsilon \BCC\rangle+\langle\BT\BP_\varepsilon \BCC,\BQ_\varepsilon \BCC\rangle\le \langle\BGL,\BQ_\varepsilon \BGL\rangle^{1/2}\langle\BCC, \BQ_\varepsilon\BCC \rangle^{1/2}\;.
\end{equation*}
Next Lemma \ref{lemest_innerprod} provides a lower bound for the left-hand side and consequently we have
\begin{equation*}\label{scalareq2}
\langle\BCC, \BQ_\varepsilon \BCC\rangle^{1/2}\Big(1-\text{const }d^{-3}\frac{\langle \BQ_\varepsilon\BCC,  \BQ_\varepsilon\BCC\rangle}{\langle\BCC, \BQ_\varepsilon \BCC\rangle}\Big)\le \langle\BGL,\BQ_\varepsilon \BGL\rangle^{1/2}\;.
\end{equation*}
As the eigenvalues of $\BCQ_\varepsilon$ are $O(\varepsilon^3)$ (see (\ref{eigest}) and  (\ref{eqDeps})), from this it is possible to derive that 
\begin{equation*}\label{scalareq3}
\langle\BCC, \BQ_\varepsilon \BCC\rangle^{1/2}(1-\text{Const }\varepsilon^3 d^{-3})\le \langle\BGL,\BQ_\varepsilon \BGL\rangle^{1/2}\;.
\end{equation*}
We now recall the constraint (\ref{eqepstod}) and it then follows that the algebraic system (\ref{system2}) is solvable and from the preceding inequality it can be determined that the estimate  (\ref{sumest1}) holds. Thus the proof of Lemma \ref{lemsolvability} is complete. \hfill $\Box$
 
\section{The energy estimate for the remainder $R_N$}\label{energymesotran}
Here, in several steps, we prove the  next lemma.
\begin{lem}
 \label{thm1_alg_f_trans2} 
Let  
\[\varepsilon < c\, d\;,\]
where $c$ is a sufficiently small absolute constant. 
Then 
remainder $R_N$, in the approximation $(\ref{introeq1})$, satisfies the energy
estimate
\begin{equation}\label{introeq22}
 \mu_O \int_{\Omega_N} |\nabla R_N|^2\, d\Bx+\sum_{1 \le j\le N}\mu_{I_j}\int_{\omega_\varepsilon^{(j)}} |\nabla R_N|^2\, d\Bx\le \text{\emph{const} } \Big\{ \varepsilon^{11
 }d^{-11
 } + \varepsilon^{5}d^{-3} \Big\} \| \nabla w_f
 \|^2_{L_2( \GO 
 )}  
  .
 \end{equation}
\end{lem}

\subsubsection*{The problem for $R_N$} The formal asymptotic algorithm leading to (\ref{introeq1}) was given in Section \ref{Formalgtran}, and invertibility of the system (\ref{alg_s_intro}) was proved in the previous section. Therefore, our objective here is derive the estimate (\ref{introeq22}). From (\ref{introeq1}), and the problems of section \ref{modprob_incl}, we have that $R_N$  is a solution of the problem 
\[\mu_O \Delta R_N(\Bx)=0\;, \quad \Bx \in \Omega_N\;,\]
\[\mu_{I_j} \Delta R_N(\Bx)=0\;, \quad \Bx \in \omega_\varepsilon^{(j)}, j=1, \dots, N\;,\]
with the exterior boundary conditions
\[R_N(\Bx)=-\sum_{k=1}^N \BC^{(k)}\cdot \Big\{ \BCD^{(k)}_\varepsilon(\Bx)-\BCP^{(k)}_\varepsilon\nabla_\By H(\Bx, \By)\Big|_{\By=\Oj^{(k)}}\Big\}\;, \quad \X \in \partial \Omega\;,\]
and the transmission conditions on the interfaces of small inclusions
\begin{eqnarray}
&&{R_N(\Bx)\Big|_{\partial \omega^{(j)+}_\varepsilon}=R_N(\Bx)\Big|_{\partial \omega^{(j)-}_\varepsilon}\;,}\nonumber\\\nonumber \\&&\nonumber\jump{R_N}{{j}}\\
&=&-(\mu_O-\mu_{I_j}) \N^{(j)} \cdot \Big\{ \nabla w_f(\Bx)+\BC^{(j)}-(\nabla_{\Bx} \otimes \nabla_{\By}) H(\Bx, \Oj^{(j)}) \BCP^{(j)}_\varepsilon \BC^{(j)}\nonumber\\ &&+\sum_{\substack{k \ne j\\ 1 \le k \le N}} \nabla_{\Bx} ( \BC^{(k)} \cdot \{ \BCD^{(k)}_\varepsilon(\Bx)-\BCP^{(k)}_\varepsilon \nabla_\By H(\Bx, \Oj^{(k)})\})\Big\}\;,
\label{rNtrancond}
\end{eqnarray}
for $1\le j \le N$.
The right-hand side of condition (\ref{rNtrancond}) is also satisfies
\begin{equation}\label{sbrn}
\int_{\partial \omega_\varepsilon^{(j)}}\left\{\jump{R_N}{{j}}\right\}\, dS_\Bx=0\;, \quad 1\le j\le N\;.
\end{equation}

\subsubsection*{Auxiliary functions }In order to derive the energy estimate for $R_N$, we need to construct functions $\Psi_k$, $k=0, \dots, N,$ such that the conditions 
\begin{equation}\label{enesttran1}
R_N(\Bx)+\Psi_0(\Bx)=0, \quad \Bx \in \partial \Omega \quad \text{ and }  
\end{equation}
\begin{equation}\label{enesttran2}
\mu_O \Big[\Dn{R_N}{}(\Bx)+\Dn{\Psi_j}{}(\Bx)\Big]\Big|_{\partial \omega_\varepsilon^{(j)+}}-\mu_{I_j} \Big[\Dn{R_N}{}(\Bx)+\Dn{\Psi_j}{}(\Bx)\Big]\Big|_{\partial \omega_\varepsilon^{(j)-}}=0\;, \quad \text{ for }j=1, \dots, N\;.
\end{equation}
In view of the boundary condition (\ref{bcH}) for $H$, we choose $\Psi_0$ in the form
\[\Psi_0(\Bx)=\sum_{1 \le k\le N} \BC^{(k)} \cdot \Big\{\BCD^{(k)}_\varepsilon(\Bx)-\BCP^{(k)}_\varepsilon  \frac{(\Bx-\BO^{(k)})}{4\pi \mu_O|\Bx-\BO^{(k)}|^3}\Big\}\;,\]
and making use of the algebraic equations (\ref{alg_s_intro}), allows $\Psi_k$ to take the representation
\begin{eqnarray*}
\Psi_k(\Bx)&=& w_f(\Bx)-w_f(\BO^{(k)})-(\Bx-\BO^{(k)})\cdot \nabla w_f(\BO^{(k)})\\\nonumber \\
&& - \BC^{(k)}\cdot \BCP^{(k)}_\varepsilon \nabla_\By H(\Bx, \Oj^{(k)})+ \sum_{\substack{ j\ne k\\ 1 \le j \le N}} \BC^{(j)} \cdot \{ \BCD^{(j)}_\varepsilon(\Bx)-\BCP^{(j)}_\varepsilon \nabla_\By H(\Bx, \By)\Big|_{\By=\BO^{(j)}}\}\\
&&-\sum_{\substack{j \ne k\\ 1 \le j \le N}} (\Bx-\BO^{(k)}) \cdot (\nabla_{\Bz} \otimes \nabla_{\Bw}) G(\Bz, \Bw)\Big|_{\substack{\Bz=\BO^{(k)}\\ \Bw=\BO^{(j)}}}\BCP^{(j)}_\varepsilon \BC^{(j)}\;, \quad k=1, \dots, N\;,
\end{eqnarray*} 
where it is easily checked that the above functions satisfy their respective boundary conditions (\ref{enesttran1}) and (\ref{enesttran2}). In addition to this, we note that 
\[\Delta \Psi_0(\Bx)=0\;,\quad \Bx \in \Omega_N\;,  \]
and 
\begin{equation}\label{delpsik}
\Delta \Psi_k(\Bx)=\mu_O^{-1}f(\Bx)\;, \quad \X \in  \bigcup_{j=1}^N\omega_\varepsilon^{(j)} \cup \Omega_N\;,\end{equation}
where $f(\Bx)$ is extended by zero inside $\omega$.
Each $\Psi_k$, $k=0, \dots, N$, is continuous across the frontiers of the inclusions $\omega_\varepsilon^{(j)}$, $j=1, \dots, N$. Equations (\ref{sbrn}) and (\ref{enesttran2}) show that
\begin{equation}\label{selfbpsi}
\int_{\partial \omega_\varepsilon^{(j)}}\left\{\jump{\Psi_j}{{j}}\right\}\, dS_\Bx=0\;, \quad \text{ for }j=1, \dots, N\;,
\end{equation}
which will be used in what follows.

In addition, we introduce cut-off functions, so that we may localize integrals over the domains $ \Omega_N$ and $\omega_\varepsilon^{(j)}$, $j=1, \dots, N,$ over regions that are in the immediate vicinity of their boundaries.
Let $(1 -\chi_0) \in C^\infty_0 (\Omega)$, such that this function is zero in a neighbourhood $\CV=\{\Bx: \text{dist}(\Bx, \partial \Omega)\le 1/2, \Bx \in \Omega\}$ of $\partial \Omega$ and equal to 1 over a neighbourhood of $\overline{\omega}$. The function $\chi_k \in C^\infty_0(\Omega)$ is chosen so that it is  equal to unity on $B^{(k)}_{2 \varepsilon}$ and vanishes outside $B_{3\varepsilon}^{(k)}$, $1\le k\le  N$.

\subsubsection*{Auxiliary estimate of the energy for $R_N$} 

Here we develop the proof of an auxiliary inequality which is important in proving Lemma \ref{thm1_alg_f_trans2}.

\begin{lem}\label{lemauxineq}
The inequality
\begin{eqnarray*}
 \| \nabla R_N\|^2_{L_2(\cup_{k=1}^N \omega^{(k)}_\varepsilon\cup \Omega_N )} 
 &\le& \text{\emph{const} } \Big\{ \| \Psi_0 \|^2_{L_2(\CV)}+ \| \nabla \Psi_0 \|^2_{L_2(\CV)}+ \sum_{1 \le k \le N}  \| \nabla \Psi_k\|^2_{L_2(B_{3\varepsilon}^{(k)})}
\Big\}\;.
\end{eqnarray*}
holds.
\end{lem}

\emph{Proof. }Consider the expression
\begin{eqnarray}
&&\mu_O \int_{\Omega_N} \nabla(R_N +\chi_0 \Psi_0)\cdot \nabla\big( R_N +\sum_{1 \le k\le N} \chi_k \Psi_k\big)\, d\Bx \nonumber \\
&&+\sum_{1 \le j\le N}  \mu_{I_j} \int_{\omega_\varepsilon^{(j)}} \nabla (R_N+\chi_0 \Psi_0)\cdot \nabla \big(R_N +\sum_{1 \le k \le N} \chi_k \Psi_k\big)\, d\Bx\;.\nonumber
\end{eqnarray} 
Our goal is to obtain from this an estimate for the energy integral appearing in the left-hand side of (\ref{introeq22}), via a linear combination of Dirichlet integrals for $\Psi_k$, $k=0, \dots,N$.

By the definitions of $\chi_k$, $k=0, \dots, N$, the above simplifies to
\[\mu_O \int_{\Omega_N} \nabla(R_N +\chi_0 \Psi_0)\cdot \nabla\big( R_N +\sum_{1 \le k\le N} \chi_k \Psi_k\big)\, d\Bx +\sum_{1 \le j\le N}  \mu_{I_j} \int_{\omega_\varepsilon^{(j)}} \nabla R_N\cdot \nabla (R_N +\Psi_j)\, d\Bx\]
Since $R_N$ is harmonic inside $\cup_{j=1}^N \omega_\varepsilon^{(j)} \cup \Omega_N$ and $\Psi_j$ is harmonic in $\omega_\varepsilon^{(j)}$, after integration by parts we obtain
\begin{eqnarray*}
&&\mu_O \int_{\Omega_N} \nabla(R_N +\chi_0 \Psi_0)\cdot \nabla\big( R_N +\sum_{1 \le k\le N} \chi_k \Psi_k\big)\, d\Bx +\sum_{1 \le j\le N}  \mu_{I_j} \int_{\omega_\varepsilon^{(j)}} \nabla R_N\cdot \nabla (R_N +\Psi_j)\, d\Bx\\
&=&-\mu_O \sum_{1 \le k\le N}\int_{B_{3\varepsilon}^{(k)} \backslash \overline{\omega_\varepsilon^{(k)}}}(R_N+\chi_0 \Psi_0) \Delta (\chi_k \Psi_k)\, d\Bx+\mu_O \int_{\partial \Omega} (R_N+\chi_0 \Psi_0) \Dn{}{}\Big\{ R_N+\sum_{1 \le k\le N} \chi_k \Psi_k\Big\}\, d\Bx\\
&&+\sum_{1 \le j\le N} \int_{\partial \omega_\varepsilon^{(j)}}  R_N\Big\{\mu_O \Big[\Dn{R_N}{}+\Dn{\Psi_j}{}\Big]\Big|_{\partial \omega_\varepsilon^{(j)+}}-\mu_{I_j} \Big[\Dn{R_N}{}+\Dn{\Psi_j}{}\Big]\Big|_{\partial \omega_\varepsilon^{(j)-}}\Big\}\, dS_\Bx\;.
\end{eqnarray*}
 Here the boundary integrals over $\partial \Omega$ and $\partial \omega_\varepsilon^{(k)}$, $k=1, \dots, N$, vanish due to conditions (\ref{enesttran1}) and (\ref{enesttran2}),  respectively. Next, as a result of the fact that $\text{supp }\chi_0 \cap \text{supp }\chi_k=\varnothing$, for $k=1, \dots, N$, we derive 
 \begin{eqnarray}
 &&\mu_O \int_{\Omega_N} |\nabla R_N|^2\, d\Bx+\sum_{1 \le j\le N} \mu_{I_j} \int_{\omega_\varepsilon^{(j)}} |\nabla R_N|^2\, d\Bx
 \nonumber\\&&=-\mu_O \int_{\Omega_N} \nabla R_N \cdot \nabla \Big(\sum_{0 \le k\le N} \chi_k \Psi_k\Big)\, d\Bx-\sum_{1\le j\le N} \mu_{I_j} \int_{\omega_\varepsilon^{(j)}} \nabla R_N \cdot \nabla \Psi_j\, d\Bx\nonumber\\
 &&-\mu_O \sum_{1\le k\le N}\int_{B_{3\varepsilon}^{(k)} \backslash \overline{\omega_\varepsilon^{(k)}}} R_N \Delta (\chi_k \Psi_k)\, d\Bx\;.\label{enesttran3}
 \end{eqnarray}
 
 In what follows $\overline{R^{(k)}}$ denotes the mean value of $R_N$  on the set $B^{(k)}_{3\varepsilon}$, $1\le k \le N$.  Using the property that $\chi_k=1$ in $B_{2\varepsilon}^{(k)}$, $k=1, \dots, N$, the last integral in (\ref{enesttran3}) can be written as 
 \begin{eqnarray}\label{enesttran4}
 &&\mu_O \sum_{1 \le j\le N} \int_{B^{(j)}_{3\varepsilon} \backslash \bar{\omega}_\varepsilon^{(j)}} R_N \Delta (\chi_j \Psi_j)\, d\Bx\nonumber\\
 &=&\nonumber\mu_O \sum_{1 \le j\le N} \int_{B_{3\varepsilon}^{(j)} \backslash \overline{{\omega}_\varepsilon^{(j)}}} (R_N -\overline{R^{(j)}}) \Delta(\chi_j(\Psi_j-\overline{\Psi_j}))\, d\Bx+\mu_O \sum_{1 \le j\le N} \overline{R^{(j)}}\int_{B_{3\varepsilon}^{(j)}\backslash \overline{\omega_\varepsilon^{(j)}}} \Delta (\chi_j \Psi_j)\, d\Bx\\
 &&+ \mu_O \sum_{1 \le j\le N} \overline{\Psi_j}\int_{B^{(j)}_{3\varepsilon}\backslash \overline{{\omega}_\varepsilon^{(j)}}} (R_N-\overline{R^{(j)}})\Delta \chi_j\, d\Bx+
\sum_{1\le j \le N}\overline{R^{(j)}}\mu_{I_j} \int_{{{\omega}_\varepsilon^{(j)}}} \Delta \Psi_j\, d\Bx \end{eqnarray}
 where $\overline{\Psi_k}$ is the mean value of $\Psi_j$ over $B^{(k)}_{3\varepsilon}$. Here, we have added the last term appearing on the right-hand side due to (\ref{delpsik}), and we can neglect third  term since using Green's formula gives
 \begin{eqnarray*}
 &&\nonumber\mu_O\int_{B^{(j)}_{3\varepsilon}\backslash \overline{\omega_\varepsilon^{(j)}}} (R_N-\overline{R^{(j)}})\Delta \chi_j\, d\Bx
 =\mu_O\int_{B^{(j)}_{3\varepsilon}\backslash \overline{\omega_\varepsilon^{(j)}}} (R_N-\overline{R^{(j)}})\Delta \chi_j\, d\Bx+\mu_{I_j}\int_{\omega_\varepsilon^{(j)}} (R_N-\overline{R^{(j)}})\Delta \chi_j\, d\Bx
\\& =&-\int_{\partial \omega_\varepsilon^{(j)}}\Big\{\jump{R_N}{j}\Big\}\, dS_\Bx=0\;,
 \end{eqnarray*}
 with the last equality  being  a result of (\ref{sbrn}). Noting this  and returning to (\ref{enesttran4}), we apply integration by parts  together with  (\ref{selfbpsi}) to yield
 \begin{eqnarray}\label{enesttran5}\nonumber
\mu_O \sum_{1 \le j\le N} \int_{B^{(j)}_{3\varepsilon} \backslash \overline{\omega_\varepsilon^{(j)}}} R_N \Delta (\chi_j \Psi_j)\, d\Bx &=&\nonumber\mu_O \sum_{1\le j\le N} \int_{B_{3\varepsilon}^{(j)} \backslash \overline{\omega_\varepsilon^{(j)}}} (R_N -\overline{R^{(j)}}) \Delta(\chi_j(\Psi_j-\overline{\Psi_j}))\, d\Bx
 \\
 &&\nonumber+\sum_{1 \le j\le N} \overline{R^{(j)}}\int_{ \partial \omega_\varepsilon^{(j)}} \Big\{\jump{\Psi_j}{j}\Big\}\, dS_\Bx\\
 &=&\mu_O \sum_{1 \le j\le N} \int_{B_{3\varepsilon}^{(j)} \backslash \overline{\omega_\varepsilon^{(j)}}} (R_N -\overline{R^{(j)}}) \Delta(\chi_j(\Psi_j-\overline{\Psi_j}))\, d\Bx\;.
 \end{eqnarray}
 In addition, concerning the first integral in the right-hand side of (\ref{enesttran3}), this is equivalent to
 \begin{eqnarray}
\mu_O \int_{\Omega_N} \nabla R_N \cdot \nabla \Big(\sum_{0 \le k\le N} \chi_k \Psi_k\Big)\, d\Bx&=& \mu_O\int_{\Omega_N} \nabla R_N \cdot \nabla \Big( \chi_0 \Psi_0\Big)\, d\Bx\nonumber\\
 &&+ \mu_O\sum_{1 \le k\le N}\int_{\Omega_N} \nabla R_N \cdot \nabla \Big( \chi_k (\Psi_k-\overline{\Psi_k})\Big)\, d\Bx\;,\label{enestspecial}
 \end{eqnarray}
 where $\overline{\Psi_k}$  is the mean value of $\Psi_k$ over the ball $B^{(k)}_{3\varepsilon}$,
 and this follows as a result of (\ref{sbrn}), the definition of $R_N$ and
 \[\mu_O\int_{\Omega_N} \nabla R_N \cdot \nabla  \chi_k d\Bx=\mu_O\int_{\partial \omega^{(k)}_\varepsilon}\frac{\partial R_N}{\partial n}\Big|_{\partial \omega^{(k)+}_\varepsilon}  dS_\Bx=\mu_{I_k}\int_{\partial \omega^{(k)}_\varepsilon}\frac{\partial R_N}{\partial n}\Big|_{\partial \omega^{(k)-}_\varepsilon}  dS_\Bx=0\;, \]
 for $1\le k\le N$.
The combination of (\ref{enesttran3}),  (\ref{enesttran5}) and (\ref{enestspecial}) provides the inequality
\begin{equation}\label{enesttran5a}\mu_O \int_{\Omega_N} |\nabla R_N|^2\, d\Bx+\sum_{1 \le j\le N} \mu_{I_j} \int_{\omega_\varepsilon^{(j)}} |\nabla R_N|^2\, d\Bx\le S_1+S_2
\end{equation}
where the terms $S_i$, $i=1, 2$, are given by
\begin{eqnarray}
S_1&=&\mu_O \Big|\int_{\CV} \nabla R_N \cdot \nabla \Big( \chi_0 \Psi_0\Big)\, d\Bx\Big| +\sum_{1 \le j\le N} \mu_{I_j} \Big|\int_{\omega_\varepsilon^{(j)}} \nabla R_N \cdot \nabla \Psi_j\, d\Bx\Big|\nonumber \\
&&+ \mu_O\sum_{1 \le k\le N}\Big|\int_{B_{3\varepsilon}\backslash \overline{\omega^{(k)}_\varepsilon}} \nabla R_N \cdot \nabla \Big( \chi_k (\Psi_k-\overline{\Psi_k})\Big)\, d\Bx\Big|\;, \nonumber
\end{eqnarray}
\begin{equation}\label{Seqs}
\begin{array}{c}
\displaystyle{S_2=\mu_O \sum_{1 \le j\le N} \Big|\int_{B_{3\varepsilon}^{(j)} \backslash \overline{\omega_\varepsilon^{(j)}}} (R_N -\overline{R^{(j)}}) \Delta(\chi_j(\Psi_j-\overline{\Psi_j}))\, d\Bx\Big|\;.}
\end{array}
\end{equation}

\subsubsection*{Estimate for the term $S_1$}
Cauchy's inequality applied to $S_1$ leads to
\begin{eqnarray}
\nonumber S_1&\le& \text{const } \Big\{ \mu_O \| \nabla R_N\|_{L_2( \CV)}\| \Psi_0\|_{L_2(\CV)}+ \mu_O \| \nabla R_N\|_{L_2( \CV)} \|\nabla \Psi_0\|_{L_2(\CV)}\\
&&\nonumber+ \mu_O \Big(\sum_{1\le k\le N} \| \nabla R_N\|^2_{L_2(B_{3\varepsilon}^{(k)}\backslash \overline{\omega_\varepsilon^{(k)}})}\Big)^{1/2} \Big(\sum_{1 \le k\le N}\| \nabla ( \chi_k (\Psi_k-\overline{\Psi_k}))\|^2_{L_2(B_{3\varepsilon}^{(k)}\backslash \overline{\omega_\varepsilon^{(k)}})}\Big)^{1/2}\\
&&+ \Big(\sum_{1 \le k\le N} \mu_{I_k}\| \nabla R_N\|^2_{L_2(\omega_\varepsilon^{(k)})}\Big)^{1/2} \Big(\sum_{1 \le k\le N}\mu_{I_k}\| \nabla  \Psi_k\|^2_{L_2(\omega_\varepsilon^{(k)})}\Big)^{1/2}\Big\}\;,\label{S1ineq}
\end{eqnarray}
and this together with
\begin{equation}\label{sumL2RN}
\Big(\sum_{1 \le k\le N} \| \nabla R_N\|^2_{L_2(B_{3\varepsilon}^{(k)}\backslash \overline{\omega_\varepsilon^{(k)}})}\Big)^{1/2} \le \text{const } \| \nabla R_N \|_{L_2(\Omega_N)}
\end{equation}
then yields the following majorant for the right-hand side of
 (\ref{S1ineq})
\begin{eqnarray*}
 &&\text{const } \Big\{ \mu_O \| \nabla R_N\|_{L_2( \CV)}\| \Psi_0\|_{L_2( \CV)}+ \mu_O \| \nabla R_N\|_{L_2( \CV)} \|\nabla \Psi_0\|_{L_2( \CV)}\nonumber\\
&&+ \mu_O \| \nabla R_N  \|_{L_2(\Omega_N)}\Big(\sum_{1 \le k\le N}\|  \nabla ( \chi_k (\Psi_k-\overline{\Psi_k}))\|^2_{L_2(B_{3\varepsilon}^{(k)}\backslash \overline{\omega_\varepsilon^{(k)}})}\Big)^{1/2}\nonumber\\
&&+ \Big(\sum_{1 \le k\le N} \mu_{I_j}\| \nabla R_N\|^2_{L_2(\omega_\varepsilon^{(k)})}\Big)^{1/2} \Big(\sum_{1 \le k \le N} \mu_{I_k}\| \nabla  \Psi_k\|^2_{L_2(\omega_\varepsilon^{(k)})}\Big)^{1/2}\Big\}\;.
\end{eqnarray*}
Thus,
\begin{eqnarray*}
\label{enesttran6}
\nonumber S_1&\le& \text{const } \Big(\mu_O \int_{\Omega_N} |\nabla R_N|^2\, d\Bx+\sum_{1\le j\le N} \mu_{I_j} \int_{\omega^{(j)}_\varepsilon} |\nabla R_N|^2\, d\Bx\Big)^{1/2}
\\&& \times \Big\{ \mu_O^{1/2} \| \Psi_0\|_{L_2( \CV)}+ \mu_O^{1/2}  \|\nabla \Psi_0\|_{L_2( \CV)}+ \Big(\sum_{1 \le k\le N}\mu_{I_k}\| \nabla  \Psi_k\|^2_{L_2(\omega_\varepsilon^{(k)})}\Big)^{1/2}\nonumber\\\nonumber \\
&&+ \Big(\sum_{1 \le k\le N}\mu_O\|  \nabla ( \chi_k (\Psi_k-\overline{\Psi_k}))\|^2_{L_2(B_{3\varepsilon}^{(k)}\backslash \overline{\omega_\varepsilon^{(k)}})}\Big)^{1/2}\Big\}\;.
\end{eqnarray*}
 By  Poincar\'e's inequality, we then have
\begin{equation}\label{Poinpsik}
\| \Psi_k -\overline{\Psi_k}\|_{L_2(B^{(k)}_{3\varepsilon})} \le \text{const }\varepsilon \| \nabla \Psi_k\|_{L_2(B_{3\varepsilon}^{(k)})}\;.
\end{equation}
and this  allows for the estimate 
\begin{eqnarray}
\label{enesttran8}
\nonumber S_1&\le& \text{const } \Big(\mu_O \int_{\Omega_N} |\nabla R_N|^2\, d\Bx+\sum_{1 \le j \le N} \mu_{I_j} \int_{\omega^{(j)}_\varepsilon} |\nabla R_N|^2\, d\Bx\Big)^{1/2}
\\&& \times \Big\{ \mu_O^{1/2} \| \Psi_0\|_{L_2( \CV)}+ \mu_O^{1/2}  \|\nabla \Psi_0\|_{L_2( \CV)}\nonumber\\
&&+ \Big(\sum_{1 \le k\le N}\mu_O\| \nabla  \Psi_k\|^2_{L_2(B_{3\varepsilon}^{(k)})}\Big)^{1/2}+ \Big(\sum_{1\le k\le N}\mu_{I_k}\| \nabla  \Psi_k\|^2_{L_2(\omega_\varepsilon^{(k)})}\Big)^{1/2}\Big\}\;.
\end{eqnarray}
\subsubsection*{Estimate for the term $S_2$ and proof of Lemma \ref{lemauxineq}}
Next, we return  to the term $S^{(2)}$ in (\ref{Seqs}). The Poincar\'e inequality (\ref{Poinpsik}) with $\Psi_k$ replaced by $R_N$, in conjunction with Minkowski's inequality leads to the estimate
\begin{eqnarray*}
S^{(2)} &\le& \mu_O \sum_{1 \le k \le N} \| R_N-\overline{R^{(k)}}\|_{L_2(B^{(k)}_{3\varepsilon})} \| \Delta(\chi_k (\Psi_k-\overline{\Psi_k}))\|_{L_2(B^{(k)}_{3\varepsilon} \backslash \bar{\omega}_\varepsilon^{(k)})}\\
&\le & \mu_O\sum_{1 \le k \le N} \varepsilon \| \nabla R_N\|_{L_2(B_{3\varepsilon}^{(k)})}\{\| (\Psi_k-\overline{\Psi_k}) \Delta \chi_k\|_{L_2(B^{(k)}_{3\varepsilon})}+2 \| \nabla \chi_k \cdot \nabla \Psi_k\|_{L_2(B^{(k)}_{3\varepsilon} )}\}\\
&\le & \text{const }\varepsilon \, \mu_O  \Big(\sum_{1 \le k \le N} \| \nabla R_N\|^2_{L_2(B^{(k)}_{3\varepsilon})}\Big)^{1/2}\\
&&\times\Big(\sum_{1 \le k \le N}\{ \| (\Psi-\overline{\Psi_k}) \Delta \chi_k\|^2_{L_2(B^{(k)}_{3 \varepsilon} )} +\| \nabla \chi_k \cdot \nabla \Psi_k\|^2_{L_2(B^{(k)}_{3 \varepsilon} )} \}\Big)^{1/2}\;.
\end{eqnarray*}
A second application of inequalities (\ref{sumL2RN}) and (\ref{Poinpsik}) then yields
\begin{equation}\label{enesttran7}
S^{(2)} \le \text{const } \Big(\mu_O \| \nabla R_N\|^2_{L_2(\Omega_N)} + \sum_{1 \le j\le N} \mu_{I_j} \| \nabla R_N\|^2_{L_2(\omega_\varepsilon^{(j)})} \Big)^{1/2} \Big( \sum_{1 \le k \le N} \mu_O \| \nabla \Psi_k\|^2_{L_2(B_{3\varepsilon}^{(k)})}\Big)^{1/2}\;.
\end{equation}
Combining (\ref{enesttran8}) and (\ref{enesttran7}) in (\ref{enesttran5a}) proves Lemma \ref{lemauxineq}. \hfill $\Box$

\subsubsection*{Completion of the proof of Lemma \ref{thm1_alg_f_trans2}}
The right-hand side of the inequality in Lemma \ref{lemauxineq} can be further expanded to give the estimate
\begin{eqnarray}
\nonumber \int_{\cup_{k=1}^N\omega^{(k)}_\varepsilon \cup \Omega_N} |\nabla R_N|^2\, d\Bx
&\le& \text{const } \{  \| \Psi_0 \|^2_{L_2( \CV)}+ \| \nabla \Psi_0 \|^2_{L_2( \CV)}+ \sum_{1\le j\le 3}\CT^{(j)}+ \sum_{1\le j\le 3}\CU^{(j)}
\}\;,
\label{enesttran9end}
\end{eqnarray}
where 
\begin{eqnarray*}
&&\CT^{(1)}=\sum_{1 \le j\le N} \int_{B^{(j)}_{3\varepsilon} \backslash \overline{\omega_\varepsilon^{(j)}}} |\nabla w_f(\Bx)-\nabla w_f (\Oj^{(j)})|^2\, d\Bx\;,\\
&&\CT^{(2)}=\sum_{1 \le j\le N} \int_{B^{(j)}_{3\varepsilon} \backslash \overline{{\omega}_\varepsilon^{(j)}}} \big| \sum_{\substack{k \ne j\\ 1 \le k \le N}} \{\nabla (\BC^{(k)} \cdot \{ \BCD^{(k)} (\Bx)-\BCP^{(k)} \nabla_\By H(\Bx, \By)\Big|_{\By=\BO^{(k)}}\})\\&&\qquad\qquad\qquad\qquad\qquad-\sum_{\substack{k \ne j\\ 1 \le k \le N}}(\nabla_{\Bz} \otimes \nabla_{\Bw}) G(\Bz, \Bw)\Big|_{\substack{\Bz=\BO^{(j)}\\\Bw=\BO^{(k)}}} \BCP^{(k)} \BC^{(k)}\big|^2\, d\Bx\\
\text{ and }&&\CT^{(3)}=\sum_{1 \le j\le N} \int_{B_{3\varepsilon}^{(j)} \backslash \overline{\omega_\varepsilon^{(j)}}} | \nabla (\BC^{(j)}\cdot \BCP^{(j)} \nabla_{\By} H(\Bx, \By)\Big|_{\By=\BO^{(j)}})|^2\, d\Bx\;.
\end{eqnarray*}
Here in (\ref{enesttran9end}),  $\CU^{(j)}$, $j=1, 2, 3$ are given by $\CT^{(j)}$, $j=1, 2, 3$, with the domains of integration $B^{(k)}_{3\varepsilon} \backslash \bar{\omega}^{(k)}_\varepsilon$ replaced $\omega_\varepsilon^{(k)}$, $j=1, \dots, N$.

\subsubsection*{Estimates for $\CT^{(j)}$, $1\le j \le 3$}
We first estimate the terms which are concentrated in the vicinity of the inclusions. Taylors expansion, shows that the term $\CT^{(1)}$
does not exceed
\[\CT^{(1)} \le \text{const } \varepsilon^5 d^{-3} \max_{\substack{\X \in \bar{\omega}\\ 1 \le i, j \le 3}}\Big| \frac{\partial^2 w_f}{\partial x_i\partial x_j}\Big|^2\]
and harmonicity of $w_f$ in a neighbourhood of $\bar{\omega}$, allows one to use the local regularity result for harmonic functions \cite{GilTrud} to obtain 
\begin{equation}\label{enesttran9a}
\CT^{(1)} \le \text{const } \varepsilon^5 d^{-3} \| \nabla w_f\|_{L_2(\Omega)}^2\;.
\end{equation}
The asymptotics of the dipole fields (see (\ref{asydipeq})) at infinity leads to 
\begin{eqnarray*}
&& \sum_{\substack{ k \ne j\\ 1 \le k \le N}} \nabla ( \BC^{(k)} \cdot \{\BCD^{(k)}(\Bx)-\BCP^{(k)} \nabla_{\By} H(\Bx, \By)\Big|_{\By=\BO^{(k)}}\})\\
&=&\sum_{\substack{k \ne j\\ 1 \le k \le N}} (\nabla_{\Bz} \otimes \nabla_{\Bw}) G(\Bz, \Bw)\Big|_{\substack{\Bz=\Bx\\\Bw=\BO^{(k)}}} \BCP^{(k)} \BC^{(k)}+O\Big(\sum_{\substack{k \ne j\\ 1 \le k \le N}} \frac{\varepsilon^4 |\BC^{(k)}|}{|\Bx-\BO^{(k)}|^4}\Big)\;.
\end{eqnarray*}
This, along with Taylor's expansion of the second order derivatives of $G(\Bx, \BO^{(k)})$ about $\Bx=\BO^{(j)}$, $j \ne k$, shows $\CT^{(2)}$ is majorized by
\begin{eqnarray}\nonumber
 && \text{const } \varepsilon^8\sum_{1 \le j\le N}\int_{B^{(j)}_{3\varepsilon}\backslash \overline{\omega_\varepsilon^{(j)}}}  \Big| \sum_{\substack{k \ne j \\ 1 \le j \le N}}  |\BC^{(k)}||\BO^{(j)}-\BO^{(k)}|^{-4}\Big|^2\\\nonumber
&\le &\text{const } \varepsilon^{11} \sum_{1 \le p\le N} |\BC^{(p)}|^2 \sum_{1 \le j\le N} \sum_{\substack{k \ne j\\ 1 \le k \le N}} |\BO^{(j)}-\BO^{(k)}|^{-8}\nonumber
\end{eqnarray}
Then, this and  Lemma \ref{eqepstod} yield the inequality
\begin{eqnarray}
\nonumber\CT^{(2)}&\le& \text{const } \frac{\varepsilon^{11}}{ d^{9}} \| \nabla w_f\|^2_{L_2(\Omega)} \sum_{1 \le j\le N} \sum_{\substack{k \ne j\\ 1 \le k \le N}} \frac{1}{|\BO^{(j)}-\BO^{(k)}|^{8}}\\
\nonumber&&\text{const } \frac{\varepsilon^{11}}{ d^{9}} \| \nabla w_f\|^2_{L_2(\Omega)} \iint_{\substack{\omega \times \omega:\\|\BX-\BY|>d}} \frac{d\BX\, d\BY}{|\BX-\BY|^8}\\
\label{enesttran10}
&\le& \text{const } \frac{\varepsilon^{11}}{d^{11}} \| \nabla w_f\|^2_{L_2(\Omega)}\;.
\end{eqnarray}
Since the dipole matrix $\BCP^{(j)}$ is $O(\varepsilon^3)$ and the derivatives of $H$ are bounded in $\omega$ we have
\begin{equation}\label{enesttran11}
\CT^{(3)} \le \text{const } \varepsilon^9 \sum_{1 \le j\le N} |\BC^{(j)}|^2
\le \text{const } \varepsilon^{9} d^{-3} \| \nabla w_f\|_{L_2(\Omega)}^2\;.
\end{equation}
\subsubsection*{Proofs of Lemma \ref{thm1_alg_f_trans2} and Theorem \ref{thm1_alg_f_trans2}}
Repeating similar steps as in the derivation of (\ref{enesttran9a})--(\ref{enesttran11}), we can write the estimates
\begin{equation}\label{enesttran12}
\left. \begin{array}{c}
\CU^{(1)} \le \text{const } \varepsilon^5 d^{-3} \| \nabla w_f\|^2_{L_2(\Omega)}\;,\\ \\
\CU^{(2)} \le \text{const } \varepsilon^{11} d^{-11} \| \nabla w_f\|^2_{L_2(\Omega)}\;,\\ \\
\CU^{(3)} \le \text{const } \varepsilon^9 d^{-3} \| \nabla w_f\|^2_{L_2(\Omega)}\;.
\end{array}\right\}
\end{equation}
Next, we estimate the terms which are concentrated near the exterior boundary $\partial \Omega$.
Owing to Lemma \ref{lemdipasymp} we have
\begin{eqnarray}
\nonumber
\|\Psi_0\|^2_{L_2( \CV)} &\le& \text{const }\varepsilon^8 \int_{ \CV} \Big|\sum_{1 \le k\le N} |\BC^{(k)}||\Bx-\BO^{(k)}|^{-3}\Big|^2\, d\Bx\\
&\le & \text{const } \varepsilon^8 \sum_{1 \le k \le N} |\BC^{(k)}|^2 \sum_{1 \le k \le N} \int_{ \CV} \frac{d \Bx}{|\Bx-\BO^{(k)}|^6}\nonumber\\
&\le & \text{const } \varepsilon^8 d^{-3} \| \nabla w_f \|^2_{L_2(\Omega)}\;. \label{enesttran13}
\end{eqnarray}
Finally, we address the second term in the right-hand side of (\ref{enesttran9end}). Similar reasoning which led to (\ref{enesttran13}) gives
\begin{eqnarray}
\nonumber
\|\Psi_0\|^2_{L_2(\CV)} 
&\le & \text{const } \varepsilon^8 \sum_{k=1}^N |\BC^{(k)}|^2 \sum_{1 \le k \le N} \int_{ \CV} \frac{d \Bx}{|\Bx-\BO^{(k)}|^8}\nonumber\\
&\le & \text{const } \varepsilon^8 d^{-3} \| \nabla w_f \|^2_{L_2(\Omega)}\;. \label{enesttran14}
\end{eqnarray}
Thus from (\ref{enesttran9a})--(\ref{enesttran14}) together with (\ref{enesttran9end}), we have
\[\mu_O \int_{\Omega_N} |\nabla R_N|^2\, d\Bx+\sum_{1 \le j\le N} \mu_{I_j} \int_{\omega^{(j)}_\varepsilon} |\nabla R_N|^2\, d\Bx \le \text{const } \{\varepsilon^{11} d^{-11}+\varepsilon^5 d^{-3}\} \| \nabla w_f\|^2_{L_2(\Omega)}\]
completing the proof of (\ref{introeq22}) and Theorem \ref{thm1_alg_f_trans}.\hfill$\square$

\section{The infinite space with a cluster of small inclusions}\label{corosec}
The theoretical results of sections \ref{intromesotran}--\ref{energymesotran} can be extended to an infinite medium containing a cloud of inclusions $\omega$. In this scenario, $\Omega=\mathbb{R}^3$ and the regular part of Green's function  $H\equiv 0$.

Here, we seek the approximation of the following boundary problem
\begin{equation}\label{mesotranprob1a_infinite}
\left.\begin{array}{c}
\displaystyle{\mu_O \Delta u_N(\Bx)=f(\Bx)\;, \quad \Bx \in \mathbb{R}^3 \backslash \cup_{ k=1}^N \overline{\omega_\varepsilon^{(k)}}\;,}
\\\\
\displaystyle{\mu_{I_j} \Delta u_N(\Bx)=0\;, \quad \Bx \in \omega_\varepsilon^{(j)}, \quad 1\le j\le  N\;,}
\\\\
\displaystyle{u_N(\Bx)\Big|_{\partial \omega^{(j)+}_\varepsilon}=u_N(\Bx)\Big|_{\partial \omega^{(j)-}_\varepsilon}\;, \quad 1\le j\le  N\;,}
\\ \\
\displaystyle{\mu_O\Dn{u_N}{}(\Bx)\Big|_{\partial \omega_\varepsilon^{(j)+}}=\mu_{I_j}\Dn{u_N}{}(\Bx)\Big|_{\partial \omega_\varepsilon^{(j)-}}\;, \quad 1\le j\le  N\;,}
\\
\\
\displaystyle{u_N(\Bx)\to \phi(\Bx)\;, \text{ as } \quad |\Bx| \to \infty\;,}
\end{array}\right\}
\end{equation}
where $f(\Bx)$ satisfies the conditions outlined in section \ref{intromesotran}, and now
\begin{equation}\label{phi-cond}
\int_{\mathbb{R}^3} f(\Bx)\, d\Bx=\int_{\mathbb{R}^3} \Delta \phi (\Bx)\, d\Bx\;.
\end{equation}

To construct the approximation for $u_N$, we require  the field $w_f$ which now satisfies equation (\ref{goveqw_f})  in $\mathbb{R}^3$ and the condition 
\[w_f(\Bx)\to \phi(\Bx)\;, \text{ as } \quad |\Bx| \to \infty\;.\] 

We have the theorem:
\vspace{0.2in}
\begin{thm}
 \label{thm_infinite} 
Let  
\[\varepsilon < c\, d\;,\]
where $c$ is a sufficiently small absolute constant. Then the solution $u_N(\Bx)$ of problem (\ref{mesotranprob1a_infinite})
admits the asymptotic representation
\begin{equation}\label{introeq1inf}
u_N(\Bx)=w_f
(\Bx)+
 \sum_{1\le k\le N} \BC^{(k)}\cdot  \BCD^{(k)}_\varepsilon(\Bx) 
+R_N(\Bx)\;,
\end{equation}
where $\BC^{(k)}=(C^{(k)}_1, C^{(k)}_2, C^{(k)}_3)^T$, $1\le k \le N$, 
satisfy the solvable linear
algebraic system 
\beq
 \nabla w_f(\BO^{(j)})+\BC^{(j)}+\sum_{\substack{k \ne j\\ 1\le k\le N}} (\nabla_{\Bz}\otimes \nabla_{\Bw})((4\pi |\Bz- \Bw|)^{-1})\Big|_{\substack{\Bz=\Oj^{(j)}\\ \Bw=\BO^{(k)}}}\BCP_\varepsilon^{(k)}\BC^{(k)}=\BO\;,\quad \text{ for  }1\le j\le  N\;.\eequ{alg_s_introinf}
The 
remainder $R_N$ satisfies the energy
estimate
\begin{equation*}\label{introeq2inf}
  \int_{\cup_{k=1}^N \omega^{(k)}_\varepsilon\cup \Omega_N} |\nabla R_N|^2\, d\Bx
 \le \text{\emph{const} } \Big\{ \varepsilon^{11
 }d^{-11
 } + \varepsilon^{5}d^{-3} \Big\} \| \nabla w_f
 \|^2_{L_2( \GO 
 )}  
  \end{equation*}
\end{thm}

 The proof of the above theorem follows closely that presented in sections \ref{intromesotran}--\ref{energymesotran} with obvious modifications. \hfill $\Box$

\section{Connection to an auxiliary homogenised problem for the cluster of inclusions}\label{connection}
In this section, we derive the auxiliary problem, which can be used to represent the coefficients appearing in the asymptotic approximation  (\ref{introeq1}) that are solutions to the algebraic system (\ref{alg_s_intro}),  in the case when a periodic cloud is contained inside a body.
We begin with a description of the geometry for a periodic cloud.

\subsection{Geometric assumptions for a periodic cluster}
We now divide the cloud $\omega$ up into many small cubes $Q_d^{(j)}=\BO^{(j)}+Q_d$, with $Q_d=\{\Bx: -d/2\le x_j \le d/2, 1\le j \le 3 \}$, where now $\BO^{(j)}$, $1\le j \le N$ are arranged periodically inside $\omega$. We assume $\varepsilon$ and $d$ satisfy the constraint (\ref{epscd}) and that for all $j$, $\omega_\varepsilon^{(j)}\subset Q_d^{(j)}$. In this case, the inclusions are constructed from the same set $F_\varepsilon$, such that $\omega^{(j)}_\varepsilon=\BO^{(j)}+F_\varepsilon$, $1 \le j \le N$. Let $\Omega \backslash \cup_{j=1}^N \overline{\omega^{(j)}_\varepsilon}$ be occupied by a material with shear modulus $\mu_O$. Each small inclusion is assumed to contain the same material, i.e. $\mu_{I_k}=\mu_I$, $1\le k \le N$. In this case, the polarization tensor for each inclusion is also identical and  $\BCP^{(k)}_\varepsilon=\BCP_\varepsilon$, $1\le k \le N$.
Here, the matrix $\CP_\varepsilon$ can be (i) negative definite if $\mu_O>\mu_I$, (ii) or positive definite $\mu_O<\mu_I$, (see Lemma \ref{lemdipasymp} of section \ref{modprob_incl}). 

We assume there exists the following limit
\begin{equation}\label{CQ}
\BCQ=\lim_{d\to 0} d^{-3}{\BCP_\varepsilon}\;,
\end{equation}
and the entries of $\BCQ$ are small.

In addition, when $N\to \infty$ ($d\to 0$ and subsequently $\varepsilon \to 0$), we assume the coefficients $\BC^{(j)}$ satisfy the relation
\begin{equation}\label{Cjlim}
\lim_{d\to 0}\BC^{(j)}=-\nabla \hat{u}(\Bx)\;,
\end{equation}
where $\hat{u}$ is the solution of the auxiliary homogenised problem within the domain $\omega \cup \Omega_\omega$, where $\Omega_\omega=\Omega\backslash \overline{\omega}$.

\subsection{Connection between algebraic system and auxiliary homogenised problem}
We take the algebraic system (\ref{alg_s_intro}) and rewrite this as
\[ \nabla w_f(\BO^{(j)})+\BC^{(j)}
+\sum_{\substack{k \ne j\\ 1\le k\le N}} (\nabla_{\Bz}\otimes \nabla_{\Bw})G(\Bz, \Bw)\Big|_{\substack{\Bz=\Oj^{(j)}\\ \Bw=\BO^{(k)}}}(d^{-3}\BCP^{(k)}_\varepsilon)\BC^{(k)} d^3=\BO\;,\quad \text{ for  }j=1, \dots, N\;.\]
In taking the limit as $N\to \infty$, so that $d\to 0$ and $\varepsilon \to 0$, the Riemann sum in the preceding equation can be replaced by an integral over $\omega$. This yields
\[\nabla w_f(\Bx)-\nabla \hat{u}(\Bx)
-\int_\omega (\nabla_{\Bz}\otimes \nabla_{\Bw})G(\Bz, \Bw)\Big|_{\substack{\Bz=\Bx\\ \Bw=\By}}\BCQ \nabla \hat{u}(\By) d\By=\BO\;,\quad \text{ for  }\Bx\in \omega\;,\]
where (\ref{CQ}) and (\ref{Cjlim}) have been implemented.
Then, we apply divergence throughout this equation, and multiply through by $\mu_O$ to obtain
\[-\mu_O\Delta \hat{u}(\Bx)+f(\Bx)+
\int_\omega \nabla_{\Bw}\delta(\Bx-\Bw)\Big|_{{ \Bw=\By}}\cdot \BCQ \nabla \hat{u}(\By) d\By=0\;,\quad \text{ for  }\Bx\in \omega\;,\]
from problems 1 and 2 of section \ref{modprob_incl}. 
Next we assume $\Bx \in \Omega_\omega \cup \omega$, and thus we retrieve the equation
\[-\mu_O\Delta \hat{u}(\Bx) +
\chi_\omega(\Bx)\nabla \cdot \BCQ \nabla \hat{u}(\Bx) +f(\Bx)=0\;,\quad \text{ for  }\Bx\in \Omega_\omega \cup \omega\;,\]
where $\chi_\omega$ is the characteristic function of $\omega$ and $f(\Bx)$ is zero in a neighborhood of $\omega$.
\subsection{Auxiliary homogenised problem}
Above, we derived the auxiliary homogenised equation for the body $\Omega$ containing a large periodic cluster of inclusions inside $\omega$. We now state the transmission problem for an  effective  medium representing the body with a periodic cluster of inclusions.

This governing equation for  $\hat{u}$ in $\Omega_\omega$ is
\begin{equation}\label{heq1}
\mu_O \Delta \hat{u}(\Bx) =f(\Bx)\;,\quad \text{ for  }\Bx\in \Omega_\omega\;,
\end{equation}
whereas in $\omega$ we have 
\begin{equation}\label{heq2}
\nabla \cdot (\mu_O \BI -
 \BCQ )\nabla \hat{u}(\Bx) =0\;,\quad \text{ for  }\Bx\in \omega\;,
\end{equation}
with $\BI$ being the $3\times 3$ identity matrix.

On the exterior $\partial \Omega$, we supply the Dirichlet condition
\begin{equation}\label{heq3}
\hat{u}(\Bx)=\phi(\Bx)\;, \quad \Bx \in \partial \Omega\;,
\end{equation}
and on the interface $\partial \omega$ we set the effective transmission conditions:
\begin{equation}\label{heq4}
\hat{u}(\Bx)\Big|_{\partial \omega^+}=\hat{u}(\Bx)\Big|_{\partial \omega^-}\;, \quad \mu_O \frac{\partial \hat{u}}{\partial n}(\Bx)\Big|_{\partial \omega^+}=\mu_O\frac{\partial \hat{u}}{\partial n}(\Bx)\Big|_{\partial \omega^-}-\Bn \cdot \BCQ \nabla\hat{u}(\Bx)\Big|_{\partial \omega^-}\;,
\end{equation}
where $\Bn$ is the unit-outward normal to $\omega$.
The matrix \begin{equation}
\label{ESM}
\mu_O \BI -
 \BCQ \;.
 \end{equation}
  appearing in (\ref{heq2}) is the effective stiiffness matrix for the periodic cluster $\omega$
Here, in general, the equation (\ref{heq2})  may describe  an anisotropic medium inside $\omega$.

The problem (\ref{heq1})--(\ref{heq4}) is useful in the case when one has a periodic cluster arranged inside $\omega$ and $N$ is large. As an alternative, one can then forego solving an $3N \times 3N$ algebraic system (\ref{alg_s_intro}) involving the unknown coefficients $\BC^{(j)}$ and make use of the approximation
\begin{equation}\label{Cgu}
\BC^{(j)}\simeq -\nabla \hat{u}(\BO^{(j)})\;, \quad 1\le j \le N\;.
\end{equation}

\subsection{Illustrative examples for  clusters with simple geometries}
\label{auxprobexp}

\subsubsection*{Effective stiffness matrix for large periodic clusters of spherical inclusions}

In the case when $F_\varepsilon$ is a sphere and contains a material with shear modulus $\mu_{I}$ of radius $\varepsilon$
the polarization tensor is diagonal and has the form
 \begin{equation*}\label{dipolematrixh}
\BCP_\varepsilon=4\pi\varepsilon^3 \mu_O\frac{\mu_{I}-\mu_O }{\mu_{I}+2\mu_O}\BI\;.
\end{equation*}
If we set $\varepsilon=b d$, with $b$ being sufficiently small, then from (\ref{CQ}), we obtain
\[\BCQ=4\pi b^3\mu_O\frac{\mu_{I}-\mu_O }{\mu_{I}+2\mu_O}\BI\;.\]
The effective stiffness matrix (see (\ref{ESM})) then takes the form $\hat{\mu}_I\BI$ where $\hat{\mu}_I$ is is the effective shear modulus of the cluster:
\begin{equation}\label{hmu}
\hat{\mu}_I=\mu_O\Big(1-4\pi b^3\frac{\mu_{I}-\mu_O }{\mu_{I}+2\mu_O}\Big)\;.
\end{equation}
Note that in this case the governing equation (\ref{heq2}) inside the cloud is the Laplace equation
and the  transmission conditions in (\ref{heq4}) become  
\begin{equation*}\label{heq4h}
\hat{u}(\Bx)\Big|_{\partial \omega^+}=\hat{u}(\Bx)\Big|_{\partial \omega^-}\;, \quad \mu_O \frac{\partial \hat{u}}{\partial n}(\Bx)\Big|_{\partial \omega^+}=\hat{\mu}_I \frac{\partial \hat{u}}{\partial n}(\Bx)\Big|_{\partial \omega^-}\;.
\end{equation*}

We investigate this case further in the numerical illustrations in the next section. Next, we demonstrate that in particular cases, one can construct  the explicit solution $\hat{u}$. 

\subsubsection*{Auxiliary homogenised problem a spherical inclusion $\omega$  in the infinite space}
The results of this section also apply to the case when $\Omega=\mathbb{R}^3$ discussed in section \ref{corosec} (see also Theorem \ref{thm_infinite} and the algebraic system (\ref{alg_s_introinf})).
The governing equations for the auxiliary homogenised problem, connected with the algebraic system (\ref{alg_s_introinf}), are then (\ref{heq1}), (\ref{heq2}) and (\ref{heq4}), which are also supplied with a condition at infinity:
\[\hat{u}(\Bx)\to  \phi(\Bx)\;, \text{ as } \quad |\Bx|\to \infty\;.\]
In addition,  $\phi(\Bx)$ satisfies the condition (\ref{phi-cond}) given in section \ref{corosec}.

Here we focus on the case when $\phi(\Bx)=\mu_O^{-1} x_1$ and we assume the domain $\omega$ is a sphere of radius $r$, with the centre at the origin. In addition we set $f(\Bx)\equiv 0$.
We consider a large periodic arrangement of spherical inclusions embedded inside this sphere, which are occupied by a material of shear modulus $\mu_I$. In this case, following the procedure of the previous section, one can consider a problem for a spherical inclusion $\omega$, occupied by a material having an effective shear modulus $\hat{\mu}_I$ given in (\ref{hmu}).

Therefore, we look for a solution $\hat{u}$ to the transmission problem
\begin{equation}\label{probuh}
\left.\begin{array}{l}
\displaystyle{\mu_O \Delta \hat{u}(\Bx)=0\;, \quad \Bx \in \mathbb{R}^3 \backslash \overline{\omega}\;,}\\ \\
\displaystyle{\hat{\mu}_{I} \Delta \hat{u}(\Bx)=0\;, \quad \Bx \in \omega,}\\ \\
\displaystyle{{\hat{u}}(\Bx)\Big|_{\partial \omega^{+}}={\hat{u}}(\Bx)\Big|_{\partial \omega^{-}}\;, \quad \mu_O\Dn{\hat{u}}{}(\Bx)\Big|_{\partial \omega^{+}}=\hat{\mu}_{I}\Dn{\hat{u}}{}(\Bx)\Big|_{\partial \omega^{-}}\;,}\\ \\
\displaystyle{\hat{u}(\Bx)=\mu_O^{-1} x_1 +O(|\Bx|^{-2})\quad \text{ as } |\Bx|\to \infty\;.}
\end{array}\right\}
\end{equation}
The solution $\hat{u}$ can then be constructed in the explicit  form as
\begin{equation}\label{hatsol1}
\hat{u}(\Bx)=\mu_O^{-1}x_1-\mathcal{D}_\omega(\Bx)\;,
\end{equation}
where function $\mathcal{D}_\omega$ is given by 
\begin{equation}\label{hatsol2}
\mathcal{D}_\omega(\Bx)=\left\{\begin{array}{ll}
\displaystyle{\frac{(\hat{\mu}_I-\mu_O) r^3}{\mu_O(\hat{\mu}_{I}+2\mu_O)}\frac{x_1}{|\Bx|^3} }
& \quad \text{ if }\quad \Bx\in \mathbb{R}^3\backslash \overline{\omega}\;, \\ \\
\displaystyle{\frac{(\hat{\mu}_{I}-\mu_O) }{\mu_O(\hat{\mu}_{I}+2\mu_O)} x_1}&\quad \text{ if }\quad \Bx\in  \omega\;.
\end{array}\right.
\end{equation}
We revisit this problem later in conjunction with the numerical simulations presented in the next section.

\section{Numerical illustrations}\label{numericalsimulations}
Here, we produce numerical computations that illustrate the effectiveness of the asymptotic approach investigated in this article.
We consider the case of a spherical body, containing a cluster of small spherical inclusions described in section \ref{secnumsetup}. The model solutions used in the asymptotic formula of Theorem \ref{thm1_alg_f_trans} are easily constructed in this case and they are also presented in section \ref{secnumsetup}. In particular, this consequently allows us to compare our results with benchmark finite element computations as discussed in section \ref{secnumsim}. Further, we end this section by comparing the asymptotic approximation with the solution of the auxiliary homogenised problem in section \ref{sec_homogenised_ex}, for a large periodic cluster of small spherical inclusions.

\subsection{Computational geometry and model solutions}\label{secnumsetup}
Let  $\Omega$ be the sphere of radius $R$ with centre at the origin and $\omega^{(j)}_\varepsilon$ be small spherical inclusions having centre at $\BO^{(j)}$ and radius $a^{(j)}_\varepsilon$.

In this case, $\Omega_N \cup \omega$ is a spherical body containing a region $\omega$ with spherical inclusions.
In section \ref{secnumsim}, $\omega$ is considered to be a cube and for this particular configuration we define
\begin{equation}\label{eqdef}
\varepsilon= R^{-1}\max_{k}a_\varepsilon^{(k)}\quad \text{ and }\quad d=R^{-1}\min_{1\le i,j\le N}\text{dist}(\BO^{(j)}, \BO^{(i)})\;. 
\end{equation}

In what follows, we introduce the model solutions to problems discussed in section \ref{modprob_incl}, for spherical geometries.
We begin by describing fields associated with the sphere $\Omega$.
\subsubsection*{Solutions to model problems in $\Omega$}
\emph{Solution $w_f$ to the unperturbed problem.} Let $w_f$ be a solution of the problem 1 in section \ref{modprob_incl},  where $\phi(\Bx)=0$ and $f(\Bx)$ is taken as 
\begin{equation}
f(\Bx)=\left\{ \begin{array}{ll}
R-r_f & \quad \text{ if } |\Bx|<r_f\\
0 & \quad \text{otherwise}\, . 
\end{array}\right.
\label{fsupp}
\end{equation}
Note that $\text{diam}(\text{supp }f)=2r_f$.
In this case, $w_f$ has the form:
\[w_f(\Bx)=\left\{\begin{array}{ll} 
\displaystyle{\frac{1}{6\mu_O}\Big(-\frac{|\Bx|^3}{2}+{r_f |\Bx|^2}-\frac{r_f^3(2R-r_f)}{2R}\Big)}\, ,&\quad \text{ if }|\Bx|<r_f\\ \\
\displaystyle{\frac{r_f^4}{12\mu_O}\Big(\frac{1}{R}-\frac{1}{|\Bx|}\Big)} \,,&\quad \text{otherwise.}
\end{array}\right.\]

\vspace{0.1in}\emph{Green's function in $\Omega.$} For the sphere  $\Omega$, occupied by a material with shear modulus $\mu_O$, the Green's function for the Dirichlet problem inside this set can be obtained from  (\ref{Grep}), where the regular part is given as:
\begin{equation*}
H(\Bx, \By)=\frac{1}{4\pi\mu_O} \frac{R}{|{\By}|}\frac{1}{|\Bx-\overline{\By}|}\;, \quad \text{ with } \quad \overline{\By}=\frac{R^2}{|\By|^2}\By\;.
\end{equation*}

\subsubsection*{Dipole fields for the small spherical  inclusion $\omega^{(j)}_\varepsilon$, $1\le j \le N$}
We consider a spherical inclusion $\omega^{(j)}_\varepsilon$, with centre $\BO^{(j)}$, radius $a^{(j)}_\varepsilon$ and  we assume this contains a material with shear modulus $\mu_{I_k}$. The inclusion is embedded in the infinite space which contains a material of shear modulus $\mu_O$. The vector function  $\BCD_\varepsilon^{(k)}$, whose components are the dipole fields for $\omega^{(k)}_\varepsilon$, takes the form
\begin{equation*}
\BCD_\varepsilon^{(k)}(\Bx)=\left\{\begin{array}{ll}
\displaystyle{\frac{(\mu_{I_k}-\mu_O) (a_\varepsilon^{(k)})^3}{\mu_{I_k}+2\mu_O}\frac{\Bx-\BO^{(k)}}{|\Bx-\BO^{(k)}|^3} }
& \quad \text{ if }\quad \Bx\in \mathbb{R}^3\backslash \overline{\omega^{(k)}_\varepsilon}\;, \\ \\
\displaystyle{\frac{(\mu_{I_k}-\mu_O) }{\mu_{I_k}+2\mu_O} (\Bx -\BO^{(k)})}&\quad \text{ if }\quad \Bx\in  \omega^{(k)}_\varepsilon\;.
\end{array}\right.
\end{equation*}
The polarization tensor for the small sphere is then 
 \begin{equation*}\label{PTS}
\BCP^{(k)}_\varepsilon=4\pi(a_\varepsilon^{(k)})^3\mu_O\frac{\mu_{I_k}-\mu_O }{\mu_{I_k}+2\mu_O}\BI\;.
\end{equation*}
Note, in accordance with section \ref{modprob_incl}, this matrix is negative (positive) definite when $\mu_O>\mu_{I_k}$ ($\mu_O<\mu_{I_k}$).
\begin{figure}
\subfigure[][]{ \centering
        \includegraphics[width=0.46\textwidth]{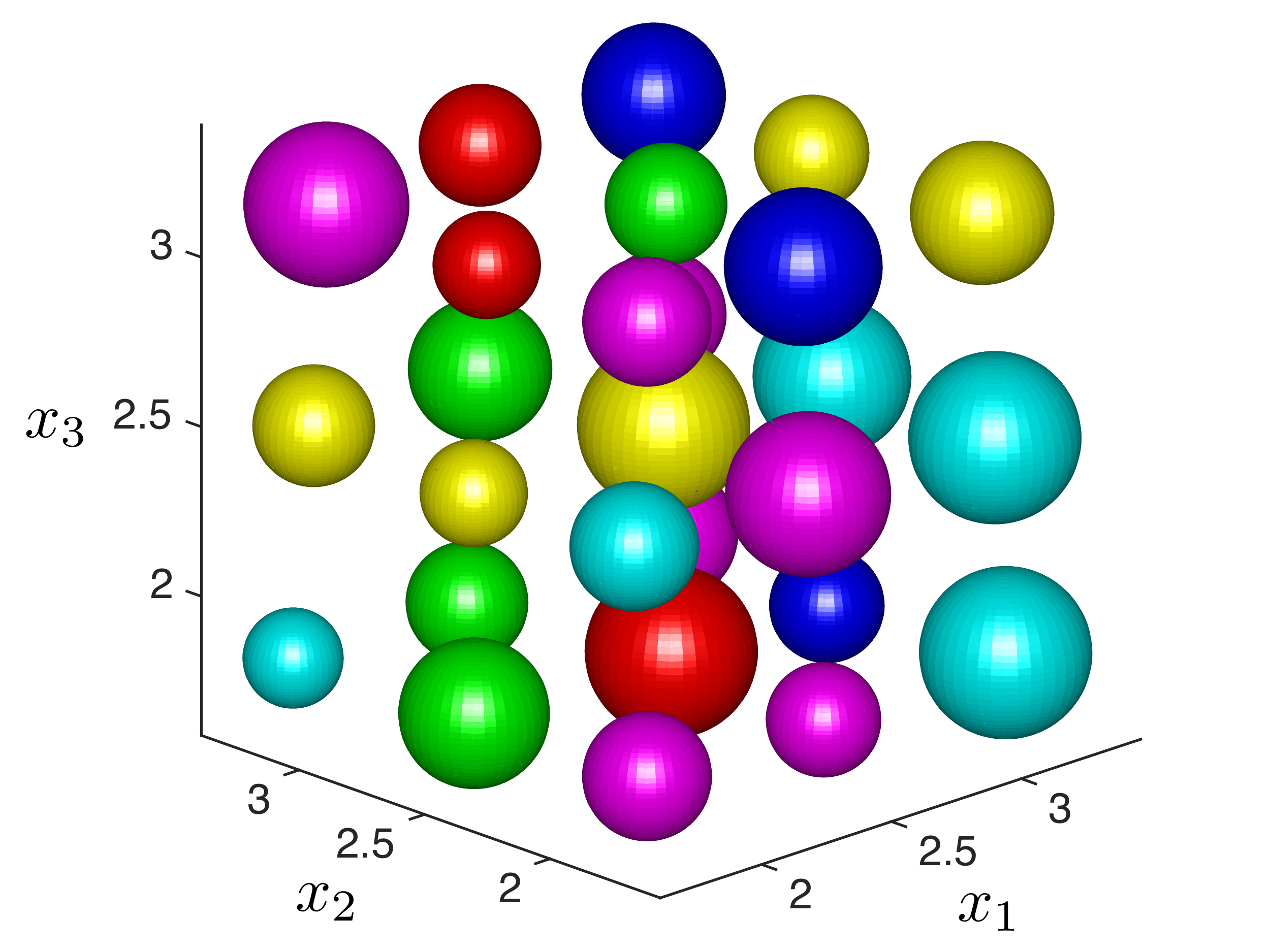}\label{fig:geom27}}
 \subfigure[][]{ \centering
        \includegraphics[width=0.46\textwidth]{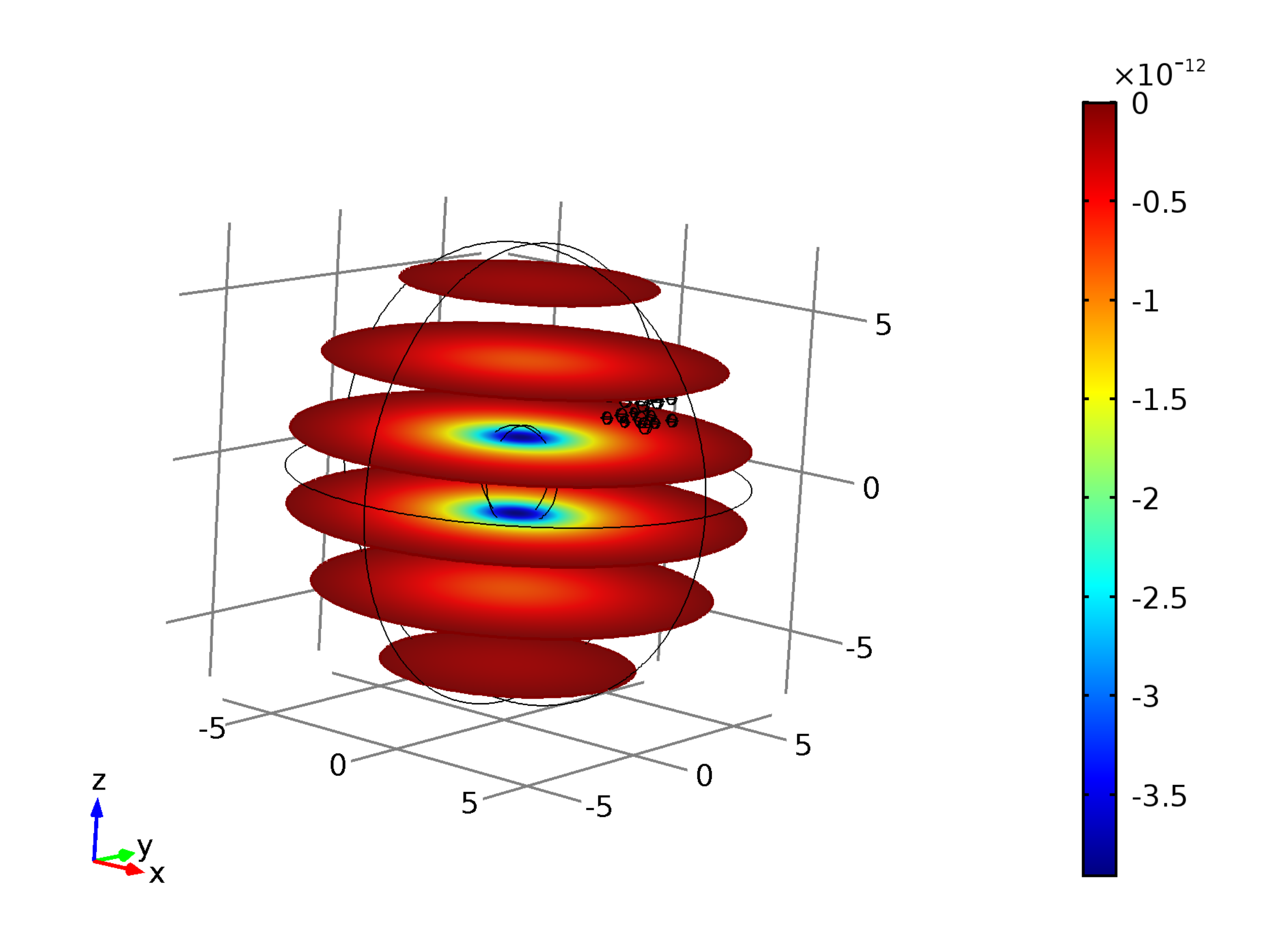}\label{fig:comsol}}      
\caption{(a) A cluster of 27 inclusions. The colors shown indicate the material contained in the inclusion: Cast Iron (green), Steel AISI 4340 (blue), Aluminum (yellow), Copper (light blue) and Iron (purple). Here, the inclusions which are red  correspond to voids (which are not occupied by a material). Properties for each inclusion can be found in Tables \ref{table:materials} and \ref{table:data}.
(b) The slice plot of the  solution of the problem outlined in sections \ref{secnumsetup} and \ref{secnumsim}, for a spherical body containing the arrangement of inclusions shown in (a). }\end{figure}

\subsection{Comparison of asymptotic approximation with the finite element method}\label{secnumsim}

For the comparison with the benchmark finite element computations in COMSOL, we consider a sphere $\Omega$ having $R=7$. The support of $f$ is contained inside the sphere of radius $r_f=1.5$. 

Inside the cluster, we assume individual inclusions are either not occupied by a material (the case when we have a void and the shear modulus inside this inclusion is set to zero) or they  contain one of the following materials: Cast Iron, Steel AISI 4340, Aluminum, Copper or Iron. The ambient matrix is occupied by Structural Steel and the  material properties used in the simulations are found in Table \ref{table:materials}.

We arrange small spherical  inclusions inside the cube $\omega$ with  centre $(2.5, 2.5, 2.5)$, having side length 2, according to the data in Table \ref{table:data}. A visual representation of the cluster, incorporating the data in Tables \ref{table:materials} and \ref{table:data},  is also shown in Figure \ref{fig:geom27}.
For this configuration of the cluster, in accordance with Table \ref{table:data} and (\ref{eqdef}), $\varepsilon=0.0343$ and $d=0.0954$.

\begin{table}[ht]
\caption{Young's modulus and Poisson's ratio of materials used in the computations.}
\centering
\begin{tabular}{|c|c|c|}
\hline\hline
Material & $\begin{array}{c}
\text{Young's Modulus},  E,\\
 (\times 10^9 \text{ N/m}^2)
\end{array}$
&Poisson's ratio, $\nu$
\\ [0.5ex] 
\hline
Cast Iron & $140$ & $0.25$\\
\hline
Steel AISI 4340 &  $205$& 0.28\\
\hline
Aluminum & $70$ & 0.33\\ 
\hline
Copper & $110$ & 0.35\\
\hline
Iron & $200 $ & 0.29\\
\hline 
Structural Steel & $200$ & 0.33\\
\hline
\end{tabular}
\label{table:materials}
\end{table}

\begin{table}[ht]
\caption{Data for inclusions contained in the cubic cloud $\omega$.}
\centering
\begin{tabular}{|c|c|c||c|c|c|}
\hline\hline
Centre & Radius & Material & Centre & Radius &  Material 
\\ [0.5ex] 
\hline
(2.56, 2.5, 1.83)& 0.24 & None &(3.12   3.15   3.16)&   0.2 & Steel AISI 4340\\
\hline
(1.83,2.48, 3.16)& 0.15& None & (1.78   2.48   2.5) & 0.15  &Aluminum \\
\hline
(2.49, 3.19, 3.16)&0.17 & None & (1.84,   3.18,   2.5)&   0.17& Aluminum\\
\hline
(1.82,  2.52, 1.83)& 0.21 & Cast Iron & (2.54,   2.51,   2.5)&   0.24 &Aluminum \\
\hline
 (2.43,   3.18,   1.83)  &  0.17& Cast Iron & (3.11,  1.83,  3.16)&   0.2 & Aluminum\\
\hline
(2.49,  3.19,  2.5)   & 0.2 & Cast Iron & (3.11,   2.51,   3.16)&   0.16 & Aluminum\\
\hline
(2.53, 2.49,   3.16)    & 0.17 & Cast Iron & (1.77,  3.19,   1.83)&   0.14  & Copper\\
\hline
(3.13  2.47   1.83)&   0.16& Steel AISI 4340 & (3.22, 1.85,   1.83)&   0.24  &Copper \\
\hline
(2.45, 1.86, 3.16)   &  0.22 &Steel AISI 4340 & (1.8,   1.86,   2.5) &   0.18 &Copper\\
\hline
\end{tabular}
\begin{tabular}{|c|c|c|}
\hline\hline
  Centre & Radius & Material
\\ [0.5ex] 
\hline
(3.13,   1.8,   2.5) &  0.24   &Copper\\
\hline
(3.15   2.47  2.5) &   0.22 &Copper\\
\hline
(1.81, 1.82,   1.83) &   0.18 &Iron\\
\hline
 (2.48, 1.81, 1.83)&   0.16 &Iron\\
\hline
(3.22,   3.19, 1.83)&  0.19 &Iron\\
\hline
 (2.44,   1.83,   2.5) &    0.23 &Iron\\
\hline
 (3.16,  3.16,   2.5)&   0.18 &Iron\\
\hline
(1.83,   1.84,   3.16) &   0.18&Iron\\
\hline
(1.85,   3.14,   3.16)&   0.23 &Iron\\
\hline
\end{tabular}

\label{table:data}
\end{table}

\begin{figure}\centering
     \subfigure[][]{
\centering
        \includegraphics[width=0.45\textwidth]{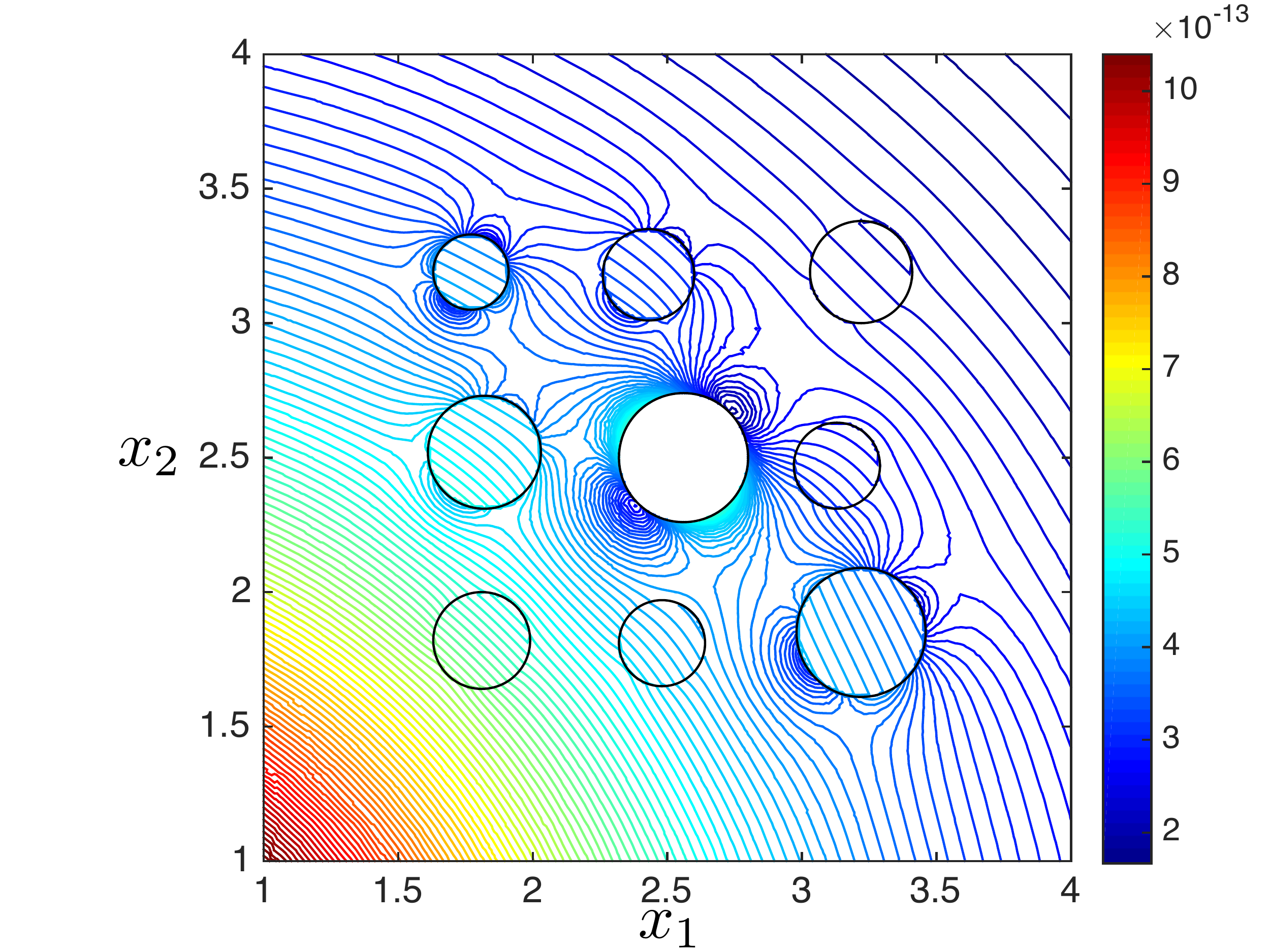}
          \label{fig:1a}}
         \subfigure[][]{
\centering
        \includegraphics[width=0.45\textwidth]{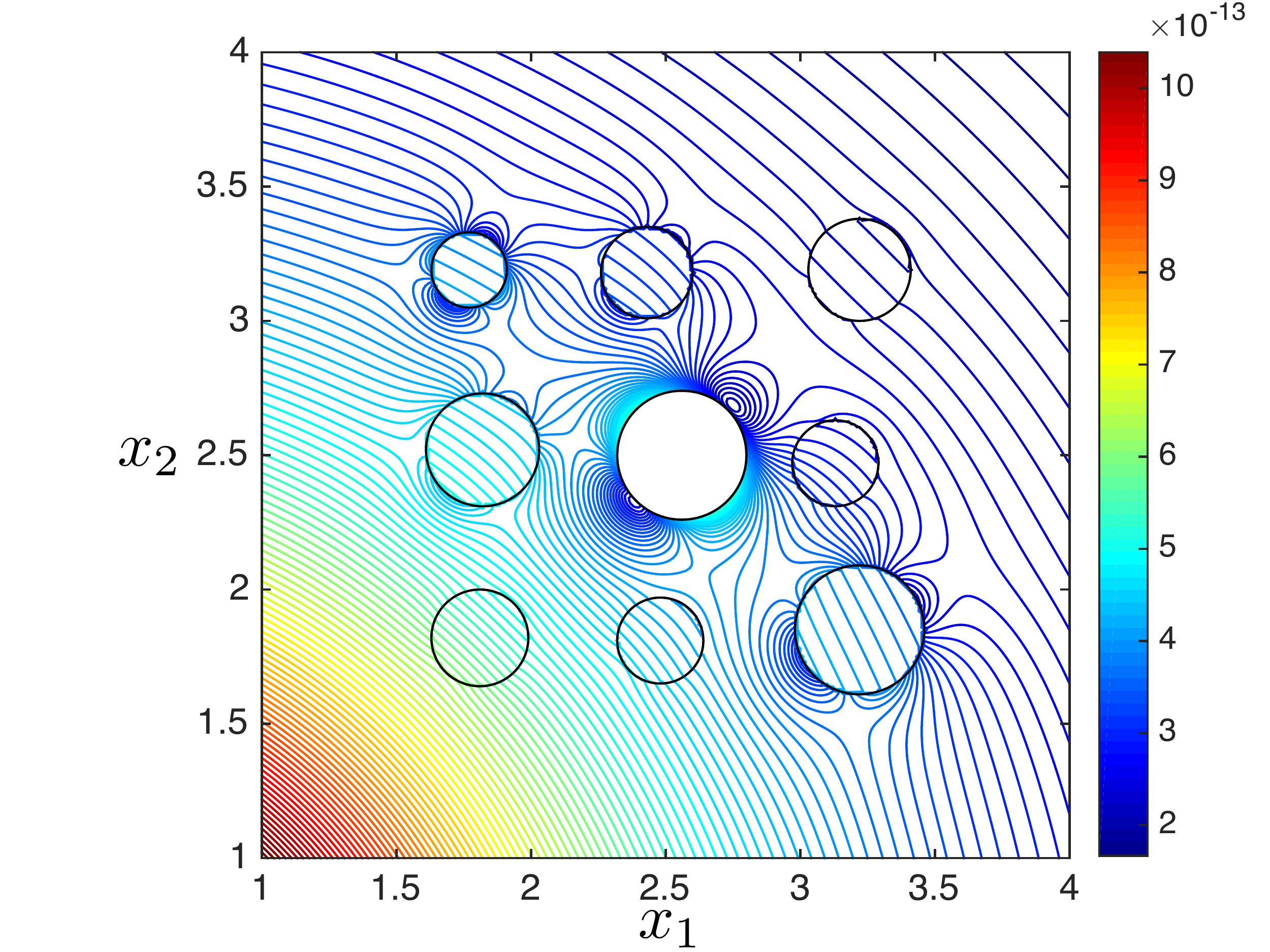}
          \label{fig:1b}}
         \subfigure[][]{
\centering
        \includegraphics[width=0.45\textwidth]{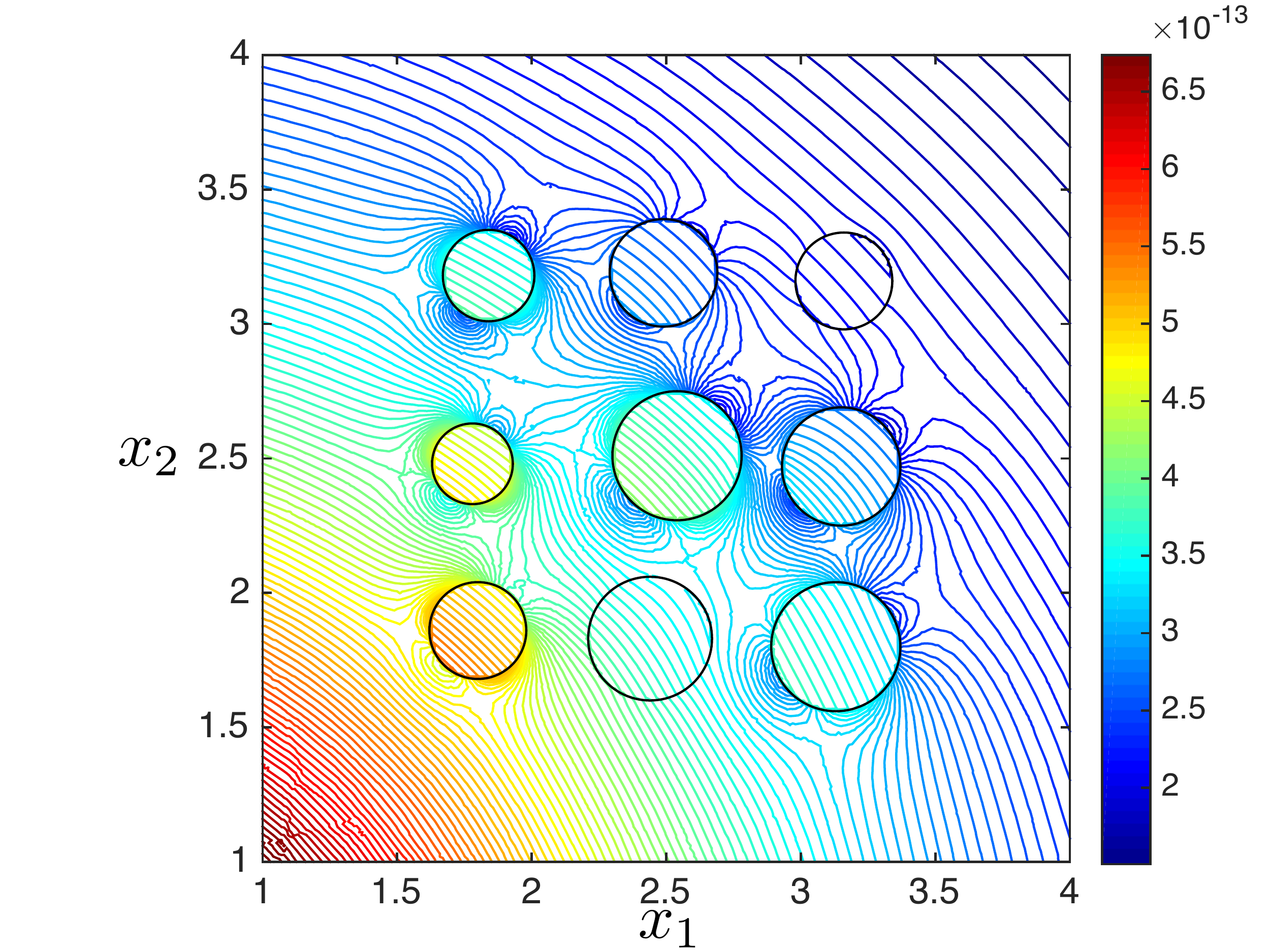}\label{fig:2a}}
           \subfigure[][]{
\centering
        \includegraphics[width=0.45\textwidth]{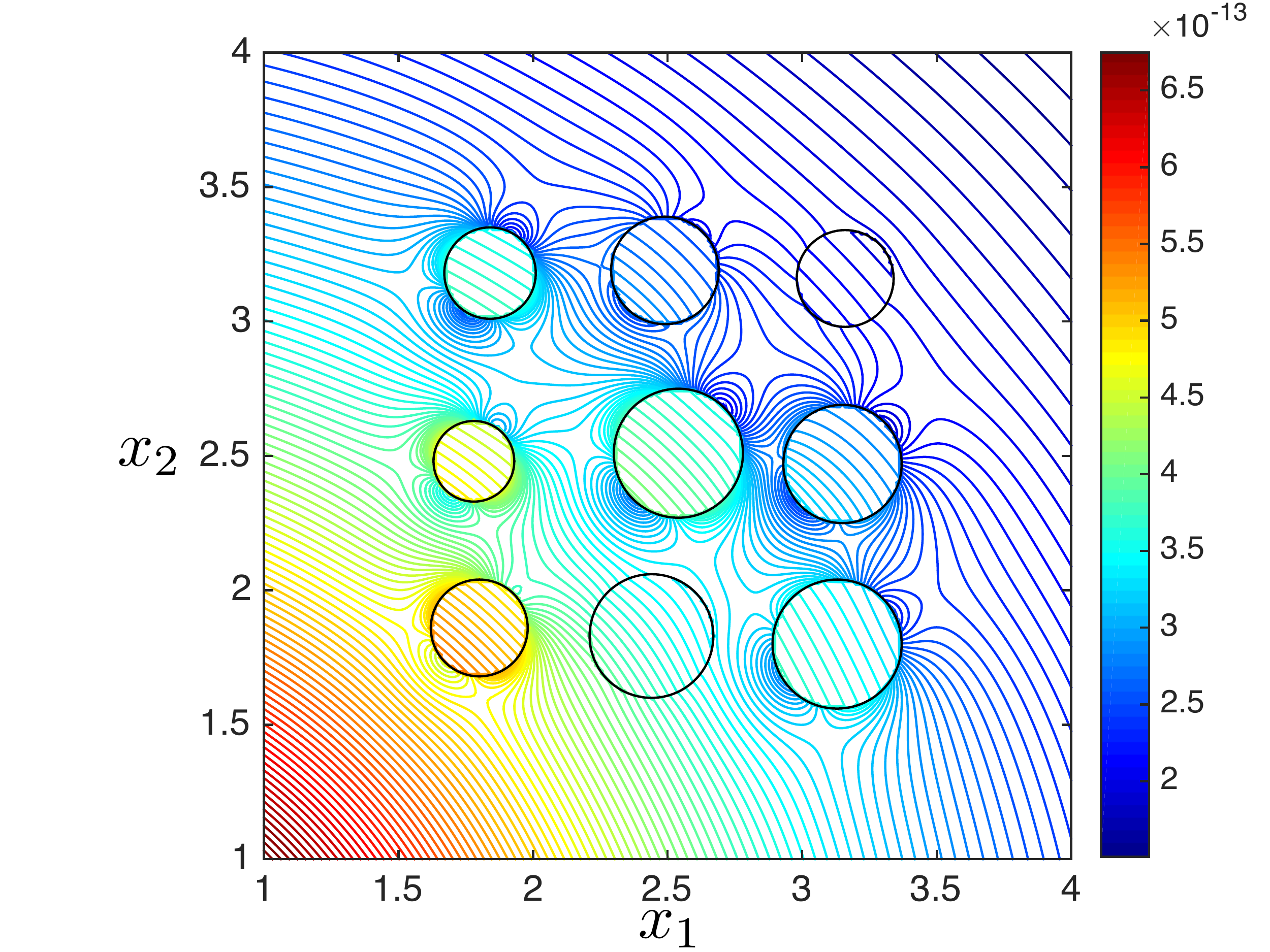}\label{fig:2b}}
        \subfigure[][]{
\centering
        \includegraphics[width=0.45\textwidth]{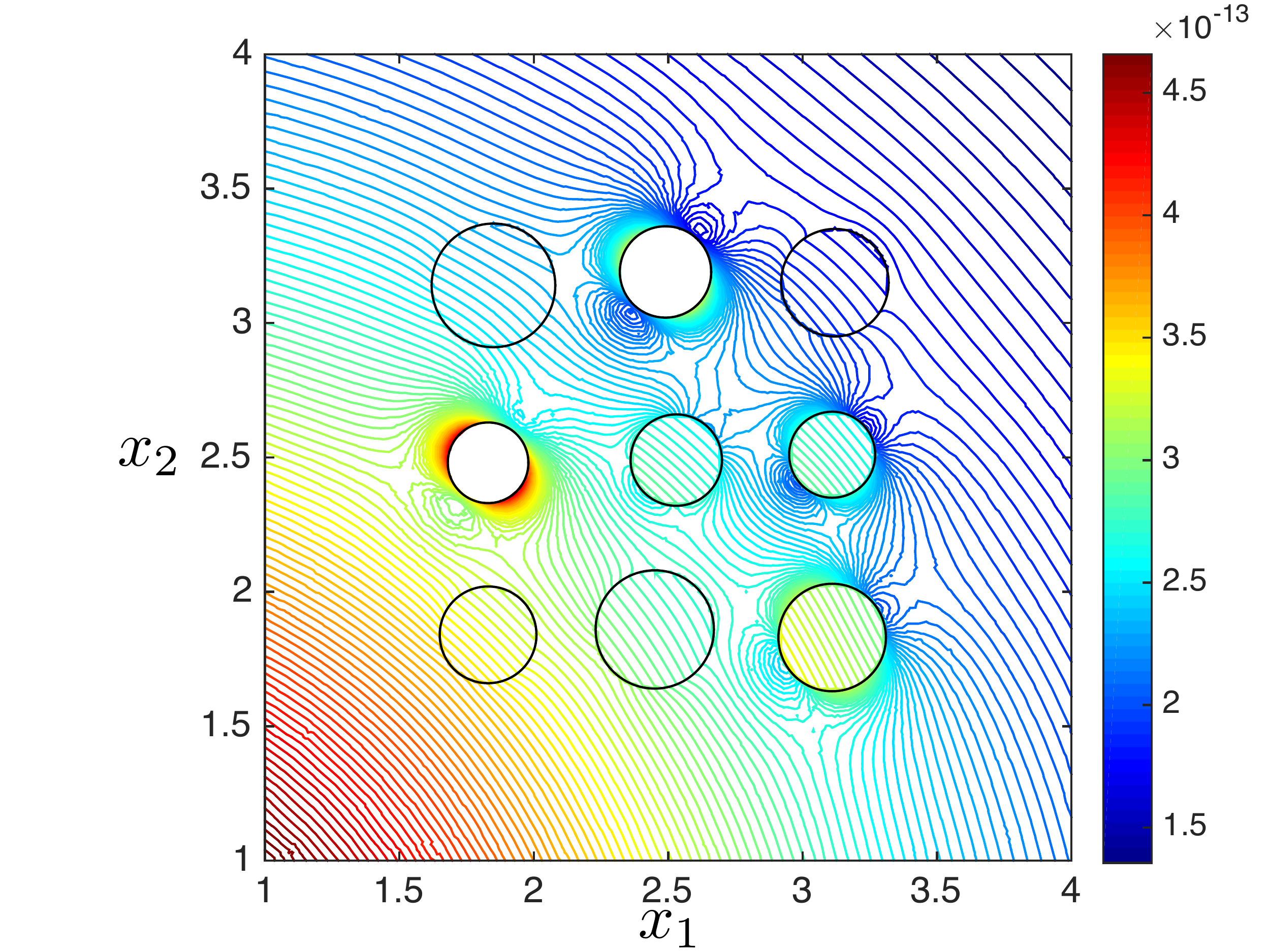}\label{fig:3a}}
           \subfigure[][]{
\centering
        \includegraphics[width=0.45\textwidth]{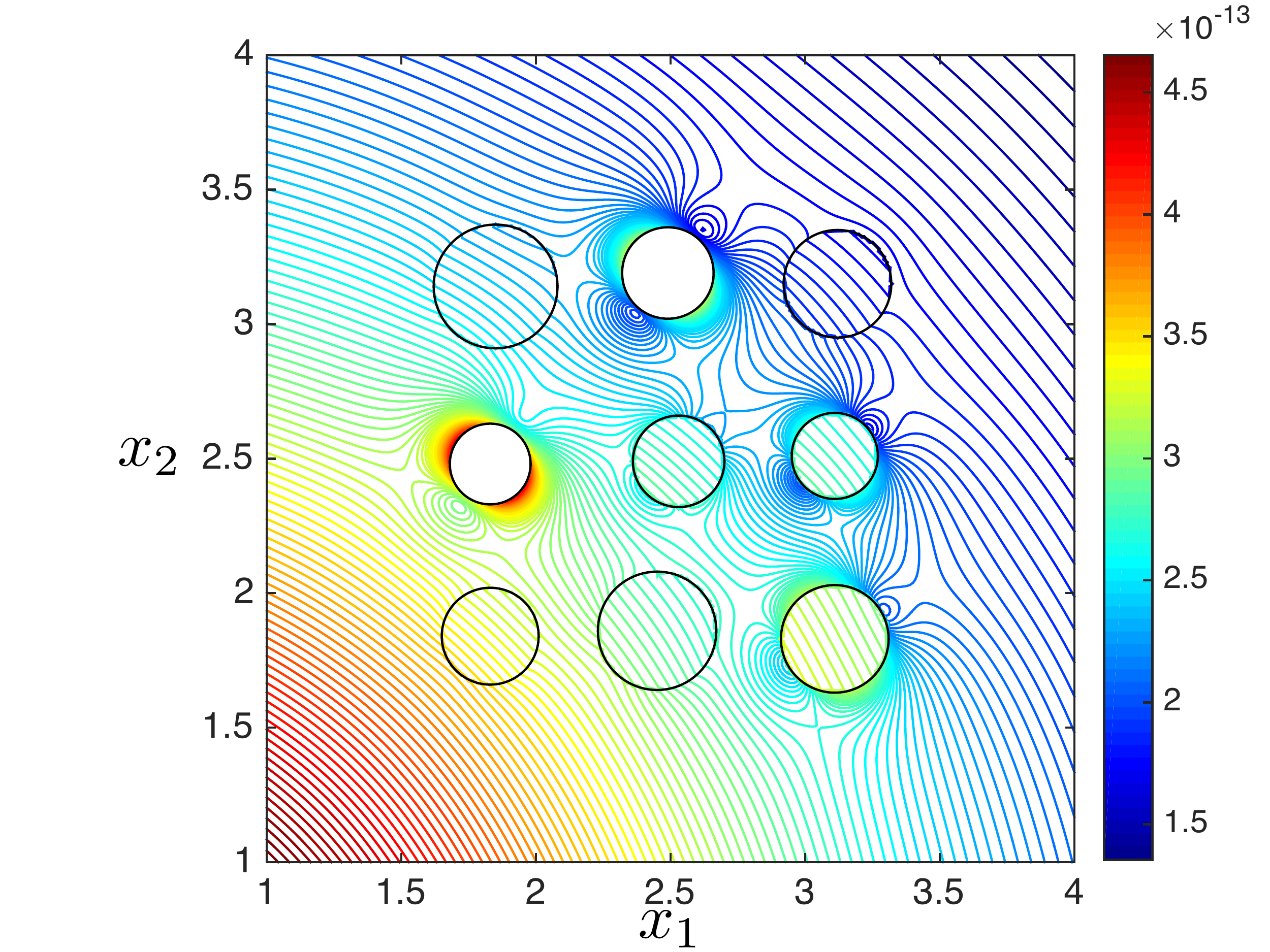}\label{fig:3b}}

\caption[]{Computations for $|\nabla u_N|$, based on COMSOL results, in the vicinity of the cluster shown in Figure \ref{fig:geom27} along cut-planes (a) $x_3=1.83$, (c) $x_3=2.5$ and (e) $x_3=3.16$. For comparison, the computations based on the asymptotic formula (\ref{introeq1}) are given in the figures on the right-hand side for (b) $x_3=1.83$, (d) $x_3=2.5$ and (f) $x_3=3.16$. }
\label{fig:2}%
\end{figure}
\subsection*{Discussion of results}
The numerical solution for $u_N$, produced by the method of finite elements in COMSOL, is shown as a slice plot in Figure \ref{fig:comsol}. Here, the effect of the non-zero support of the function $f$  (in (\ref{fsupp})) can be clearly seen inside the sphere of radius 1.5.
This computation took 1 hour 8 min and required a calculation involving 4179829 degrees of freedom.

Next we consider cut-planes that intersect the cluster and are defined by $x_3=1.83$, $x_3=2.5$ and $x_3=3.16$.
The quantity $|\nabla u_N|$ computed using the numerical solution in COMSOL,  along these cut-planes in the vicinity of the cluster, is supplied  in Figures \ref{fig:1a}  for $x_3=1.83$, \ref{fig:2a}  for $x_3=2.5$  and \ref{fig:3a} for $x_3=3.16$.   The corresponding computations for $|\nabla u_N|$ based on the derivatives of the leading order asymptotics of (\ref{introeq1}) are shown in Figures \ref{fig:1b},  \ref{fig:2b} and \ref{fig:3b}. 

 In this case, COMSOL will compute $|\nabla u_N|$ by differentiating numerically, and hence when fields are rapidly varying (for instance inside or near the cluster), one would expect some inaccuracies in the numerical results. On the other hand, formula (\ref{introeq1}) is uniform everywhere inside $\Omega_N$, and in particular, uniform up to and including the boundaries of the small inclusions. This formula can be differentiated and used to give an accurate depiction of the strain field inside the cluster.  
Here, the results  produced in COMSOL were based on a computationally intensive simulation and required an extremely fine mesh. However, we  can still find non-smooth behaviour in the strain field based on the  finite element calculations. For instance, in Figure \ref{fig:1a} at approximately $x_1=3.6$ and $x_2=1.8$ the finite element calculations appear to vary in a non-smooth fashion, indicating some slight numerical error. At this point, in Figure \ref{fig:1b}, the asymptotic formula predicts a much smoother behaviour in the strain field. Further refinement of the mesh  in COMSOL would allow one to recover the accurate behaviour of the strain field with finite elements, which would simultaneously require greater computing power.

We note that there is an excellent qualitative agreement between the computations along the cut-planes, even in this case where it is apparent the hole size is competing with their separation. In fact, the average  absolute error between the results shown in  (i) Figures \ref{fig:1a} and \ref{fig:1b} is $4.98\times 10^{-16}$,  (ii) Figures \ref{fig:2a} and \ref{fig:2b} is $5.78\times 10^{-16}$ and (iii) Figures \ref{fig:3a} and \ref{fig:3b} is $2.67\times 10^{-16}$. Thus,  there is an outstanding agreement between the results based on the numerical computations in COMSOL and those from the asymptotic approach.

Since the asymptotic formula (\ref{introeq1}) predicts the correct strain field when compared with independent finite element computations, one can use this formula for more complicated, larger clusters of spherical inclusions. For $N=64$, and the corresponding computations for $|\nabla u_N|$ are shown in Figure \ref{fig:3}. There, one would expect a more rapid variation of the strain field in a neighborhood containing the cluster. In this case, COMSOL was unable to compute the solution to this problem.

\subsection{Example: Computations for an infinite medium with a large spherical cluster of inclusions}\label{sec_homogenised_ex}

Now we consider an infinite medium containing a sphere with a periodic arrangement of many small spherical inclusions. 
We take 
$\Omega=\mathbb{R}^3$ and the domain $\omega$ as the sphere of radius $1/2$.  Inside $\omega$, we distribute $N$ small identical spherical inclusions $\omega^{(j)}_\varepsilon$. 

\subsubsection{Geometry of the spherical cluster}\label{geomspherecluster}
 We consider the cube $Q$ having side length 1 and centre at the origin. We  divide this cube into $N_1$ cubes, $Q_d^{(j)}$, described as follows.
 We introduce the set $\Sigma$ as
\[\Sigma:=\Big\{\BO_{ijk}: \BO_{ijk}=\Big(\frac{2i-1-N_1}{2N_1},\frac{2j-1-N_1}{2N_1}, \frac{2k-1-N_1}{2N_1}\Big)^T\;, 1\le i,j,k\le N_1\Big\}\;,\]
and we allow $\BP^{(j)} \in \Sigma$, $1\le j \le N_1$ such that $\text{dist}(\BP^{(j)}, \BP^{(k)})\ne 0$, for $j \ne k$, $1\le j, k\le N_1$. Setting $d=1/N_1^{1/3}$, we then have $Q_d^{(j)}=\BP^{(j)}+Q_d$, with $\BP^{(j)}$ being the centre of the cube $Q_d^{(j)}$, $1\le j \le N_1$.

\begin{figure}\centering
         \subfigure[][]{
\centering
        \includegraphics[width=0.45\textwidth]{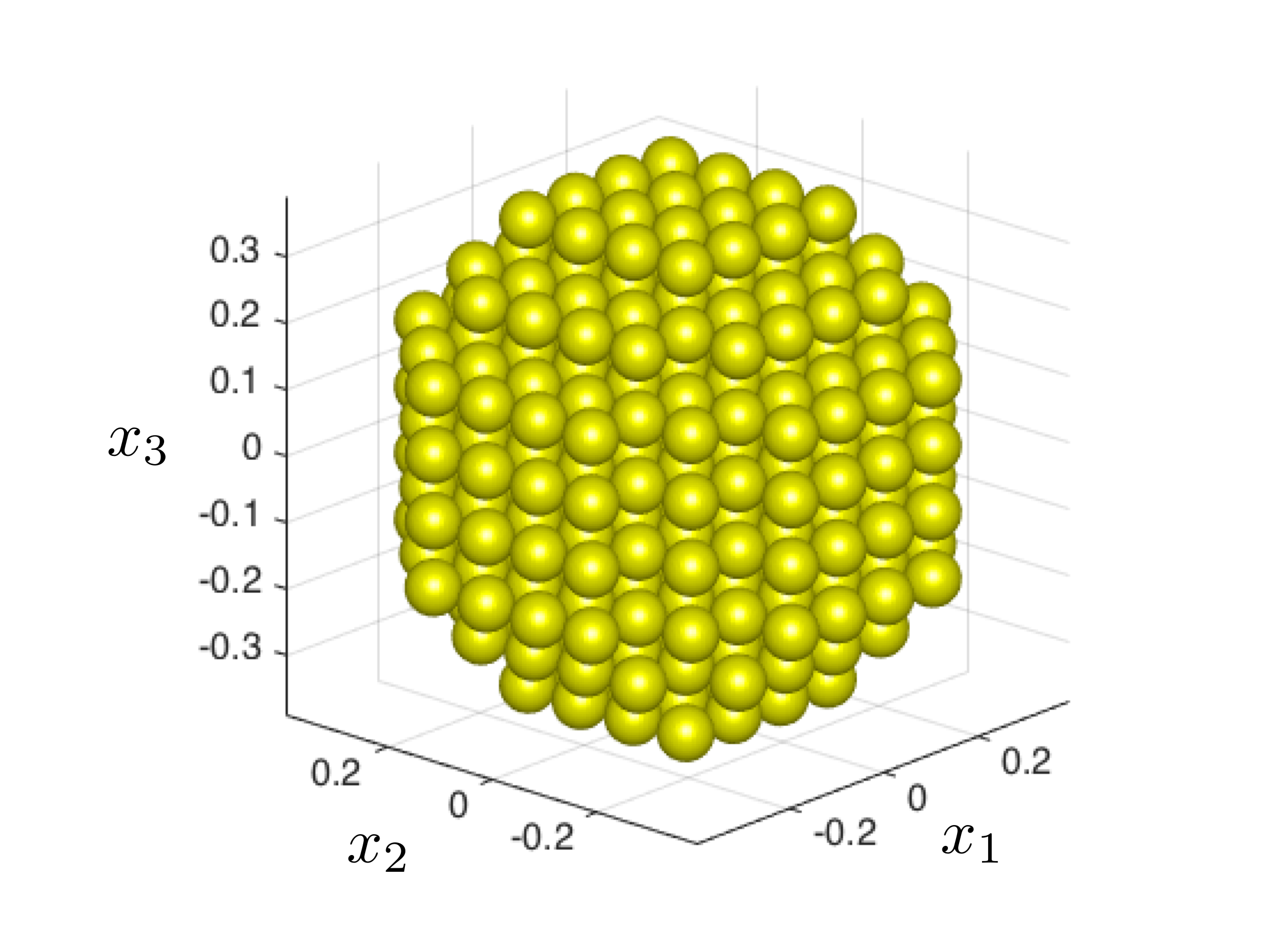}\label{fig 5a}}
           \subfigure[][]{
\centering
        \includegraphics[width=0.45\textwidth]{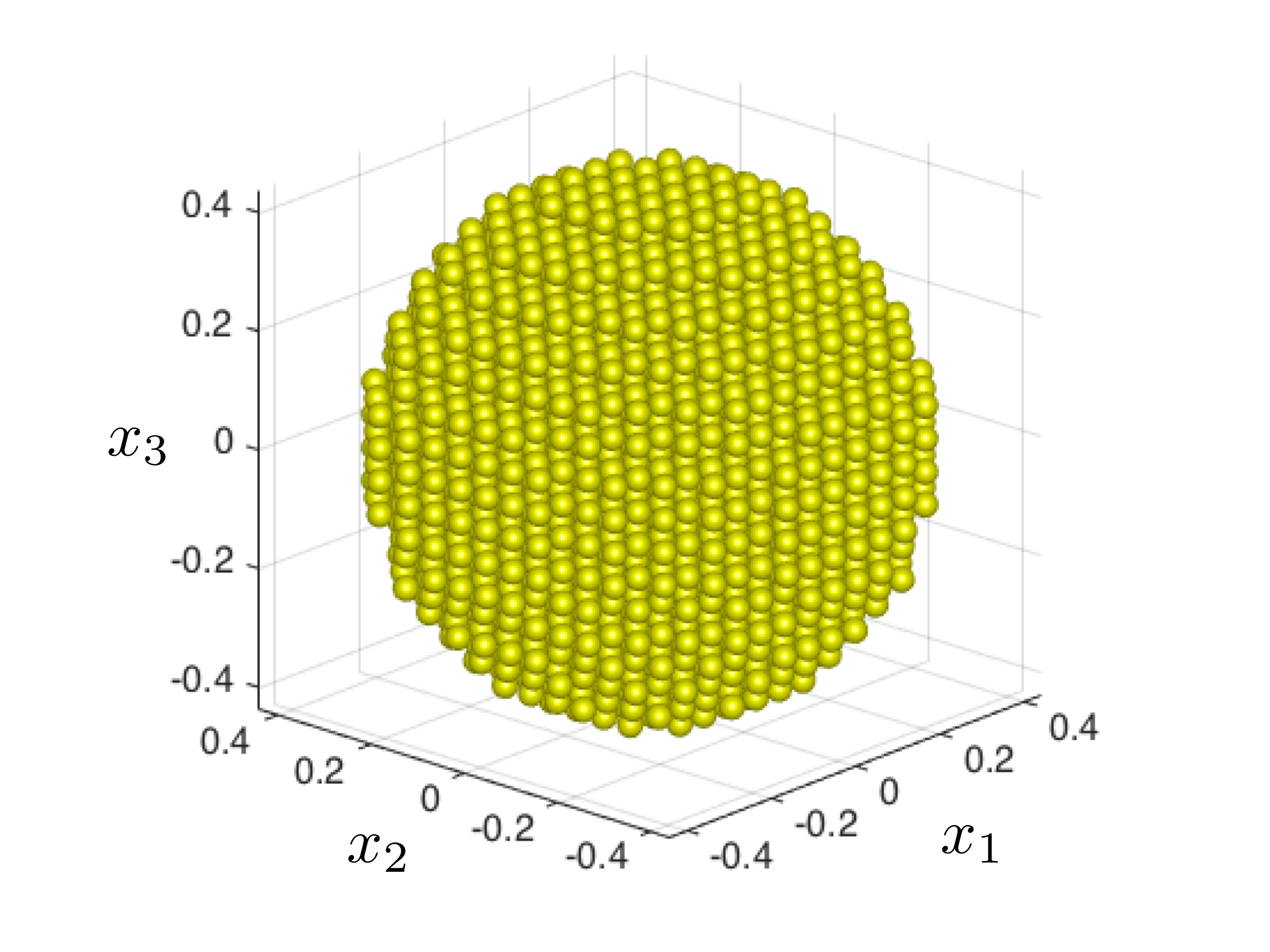}\label{fig:5b}}\\
        \subfigure[][]{
\centering
        \includegraphics[width=0.45\textwidth]{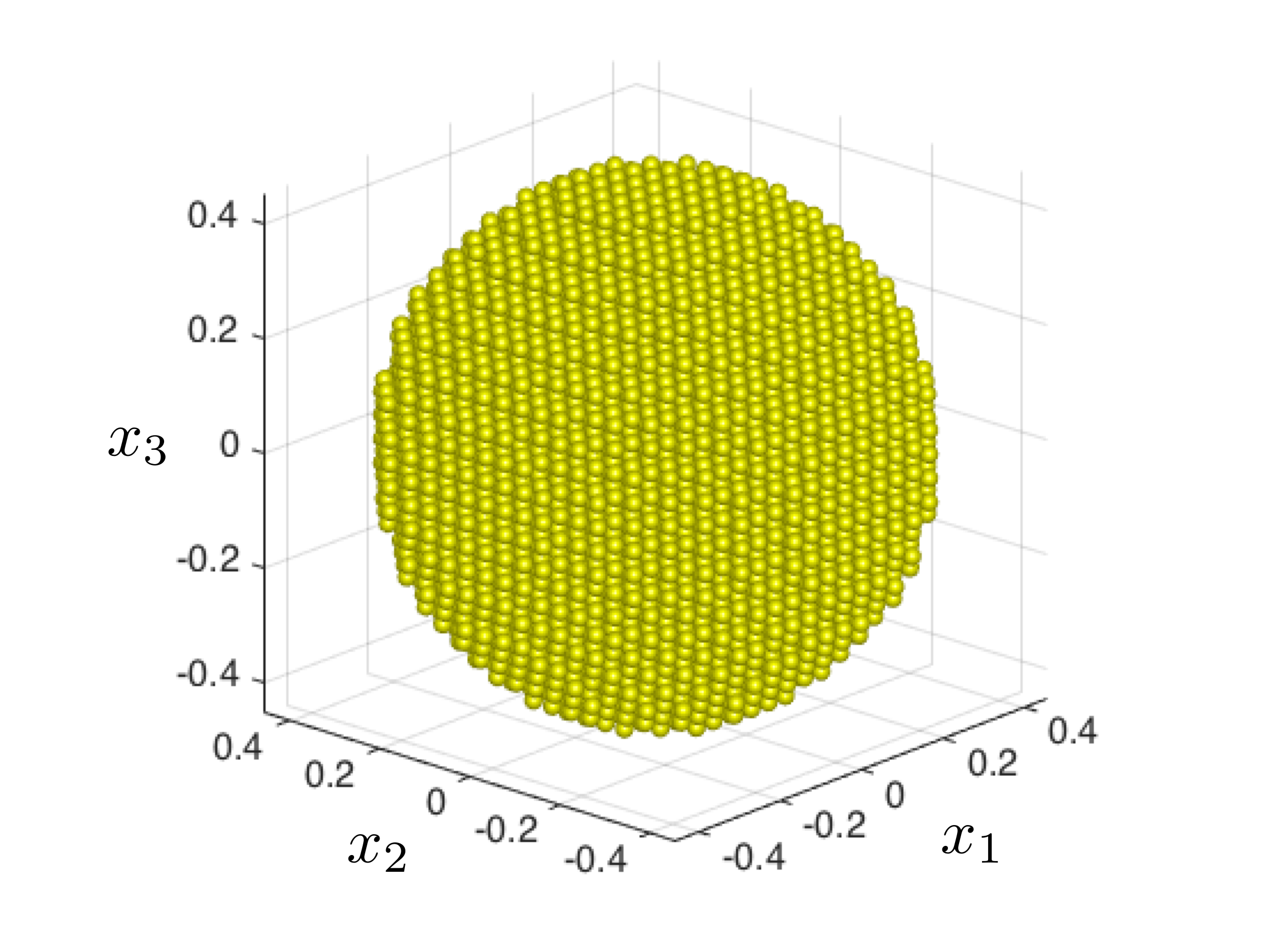}\label{fig:5c}}

\caption[]{ Clusters of spherical inclusions contained in a sphere, as described in section \ref{geomspherecluster}. Here we show the cases (a) $N=304$ ($N_1=1000$), (b) $N=2284$ ($N_1=5832$) and (c) $N=5880$, ($N_1=13824$). Each inclusion is occupied by Aluminum, whose material properties are described in Table \ref{table:materials}, in section \ref{secnumsim}.}
\label{fig:5}%
\end{figure}

To create the spherical cluster, we define the collection 
\[\Pi:=\Big\{\BP^{(j)}:Q_d^{(j)}= \BP^{(j)}+Q_d\text{ and } Q_d^{(j)}\subset \omega \text{ for } 1\le j \le N \Big\}\;.\]
we say this set has cardinality $|\Pi|=N$. Moreover, let the spherical inclusions be given by the sets 
 $\omega^{(j)}_\varepsilon=\BO^{(j)}+F_\varepsilon$, where $F_\varepsilon$ is a ball of radius $\varepsilon$ and centre at the origin. Here, $\BO^{(j)}\in \Pi$, $1\le j \le N$, such that   $\text{dist}(\BO^{(j)}, \BO^{(k)})\ne 0$, for $j \ne k$, $1\le j, k\le N$. In addition, the parameters $\varepsilon$ and $d$ are related by
\begin{equation}\label{epsdnum}
\frac{\varepsilon}{d}=b=\Big(\frac{3N_1}{4\pi N}\beta\Big)^{1/3}\;,
\end{equation}
 where $\beta<4\pi N/3N_1$ and $b$ was introduced in section \ref{auxprobexp}.
 One can verify that as $N_1 \to \infty$ that
 \[\frac{N_1}{N}\to \frac{\text{meas}(Q)}{\text{meas}(\omega)}=\frac{6}{\pi}\;,\]
 with $\text{meas}(A)$ being the three-dimensional measure of the set $A$. Typical arrangements of inclusions created according to the description provided here can be found in Figure \ref{fig:5} for $N=304$ ($N_1=1000$), $N=2284$ ($N_1=5832$)  and $N=5880$ ($N_1=13824$). In what follows, we assume the inclusions are occupied by Aluminum and the ambient matrix is occupied by Structural Steel (see Table \ref{table:materials} in section \ref{secnumsim} for the corresponding material properties).
 
\subsubsection{Governing equations for the infinite medium with a large spherical cluster}
We consider the boundary value problem:
\begin{eqnarray}\label{uninf1}
&&{\mu_O \Delta u_N(\Bx)=0\;, \quad \Bx \in \mathbb{R}^3 \backslash \cup_{ k=1}^N  \overline{\omega_\varepsilon^{(k)}}\;,}
\\
&&{\mu_{I} \Delta u_N(\Bx)=0\;, \quad \Bx \in \omega_\varepsilon^{(j)}, \quad  1\le j \le N\;,}
\\
&&{u_N(\Bx)\Big|_{\partial \omega^{(j)+}_\varepsilon}=u_N(\Bx)\Big|_{\partial \omega^{(j)-}_\varepsilon}\;, \quad 1\le j \le N\;,}
 \\
&&{\mu_O\Dn{u_N}{}(\Bx)\Big|_{\partial \omega_\varepsilon^{(j)+}}=\mu_{I}\Dn{u_N}{}(\Bx)\Big|_{\partial \omega_\varepsilon^{(j)-}}\;, \quad 1\le j \le N\;,}
\end{eqnarray}
where at infinity we prescribe 
\begin{equation}\label{uninf5}
u_N(\Bx)={\mu_O}^{-1}x_1+O(|\Bx|^{-2})\;, \text{ for } |\Bx|\to \infty\;.
\end{equation}
Here, $\mu_O$ and $\mu_I$ are the  shear moduli for Structural Steel and Aluminum, respectively.

The results of section \ref{corosec} are readily adapted to this particular boundary value problem (\ref{uninf1})--(\ref{uninf5}), by taking $w_f(\Bx)=\mu_O^{-1}x_1$. According to the algebraic system (\ref{alg_s_introinf}) 
and the procedure followed in section \ref{connection} (see (\ref{Cgu}) and section \ref{auxprobexp}), the coefficients $\BC^{(j)}$ as $N\to \infty$ (and $d\to 0$) admit the form
\begin{equation}\label{Ch}
\hat{C}=\lim_{d\to 0} C_j=-\nabla \hat{u}(\BO^{(j)})=-\mu_O^{-1}\nabla (x_1-\mathcal{D}_\omega(\Bx))\Big|_{\Bx=\BO^{(j)}}\;,
\end{equation}
where $\hat{u}(\Bx)=x_1-\mathcal{D}_\omega(\Bx)$ is the solution to the auxiliary homogenised problem stated in section \ref{auxprobexp}.

\subsubsection{Numerical comparison of asymptotic approximation with the solution to the auxiliary homogenised problem}

\begin{figure}\centering
         \subfigure[][]{
\centering
        \includegraphics[width=0.49\textwidth]{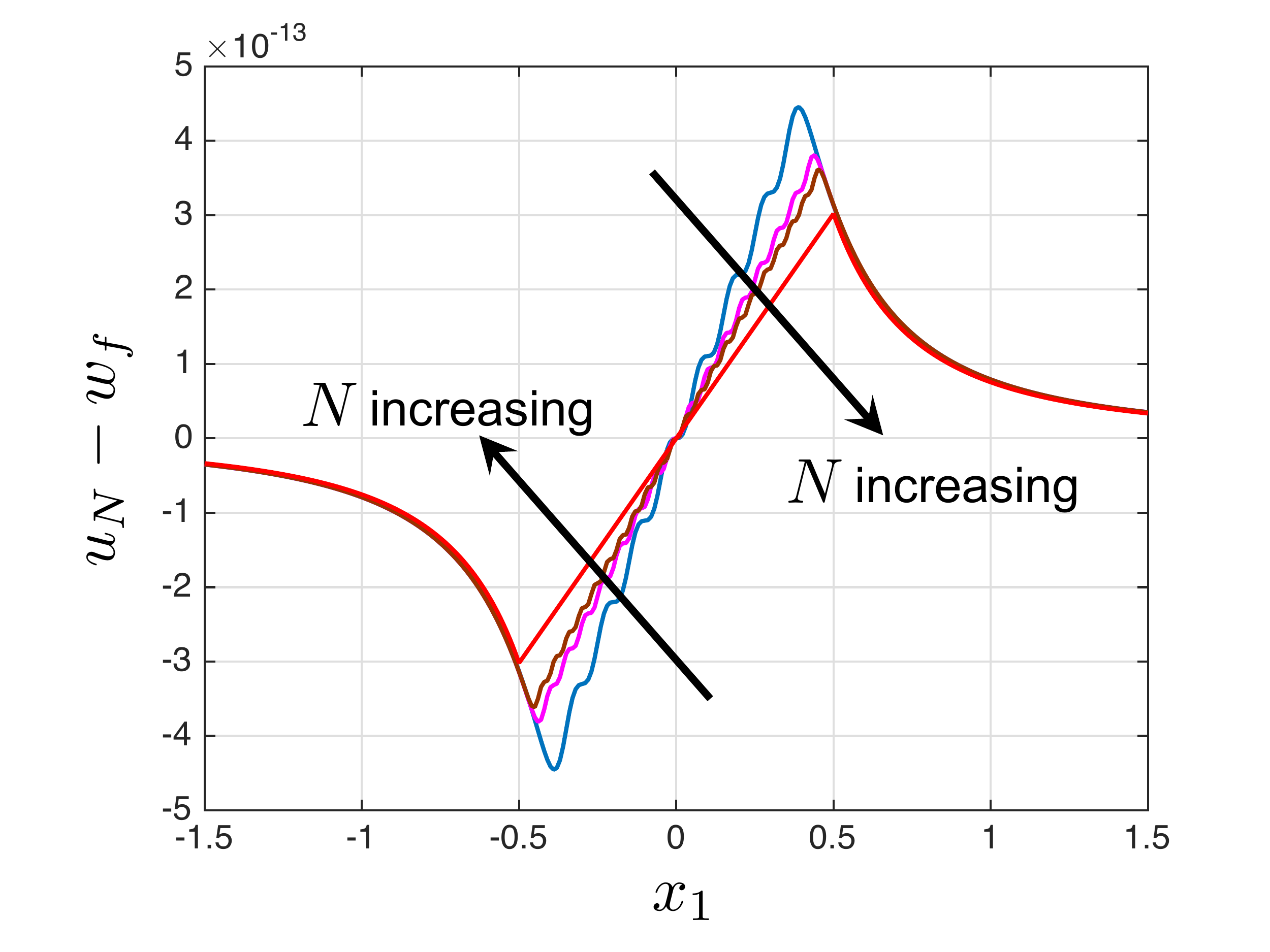}\label{fig:6a}}
           \subfigure[][]{
\centering
        \includegraphics[width=0.49\textwidth]{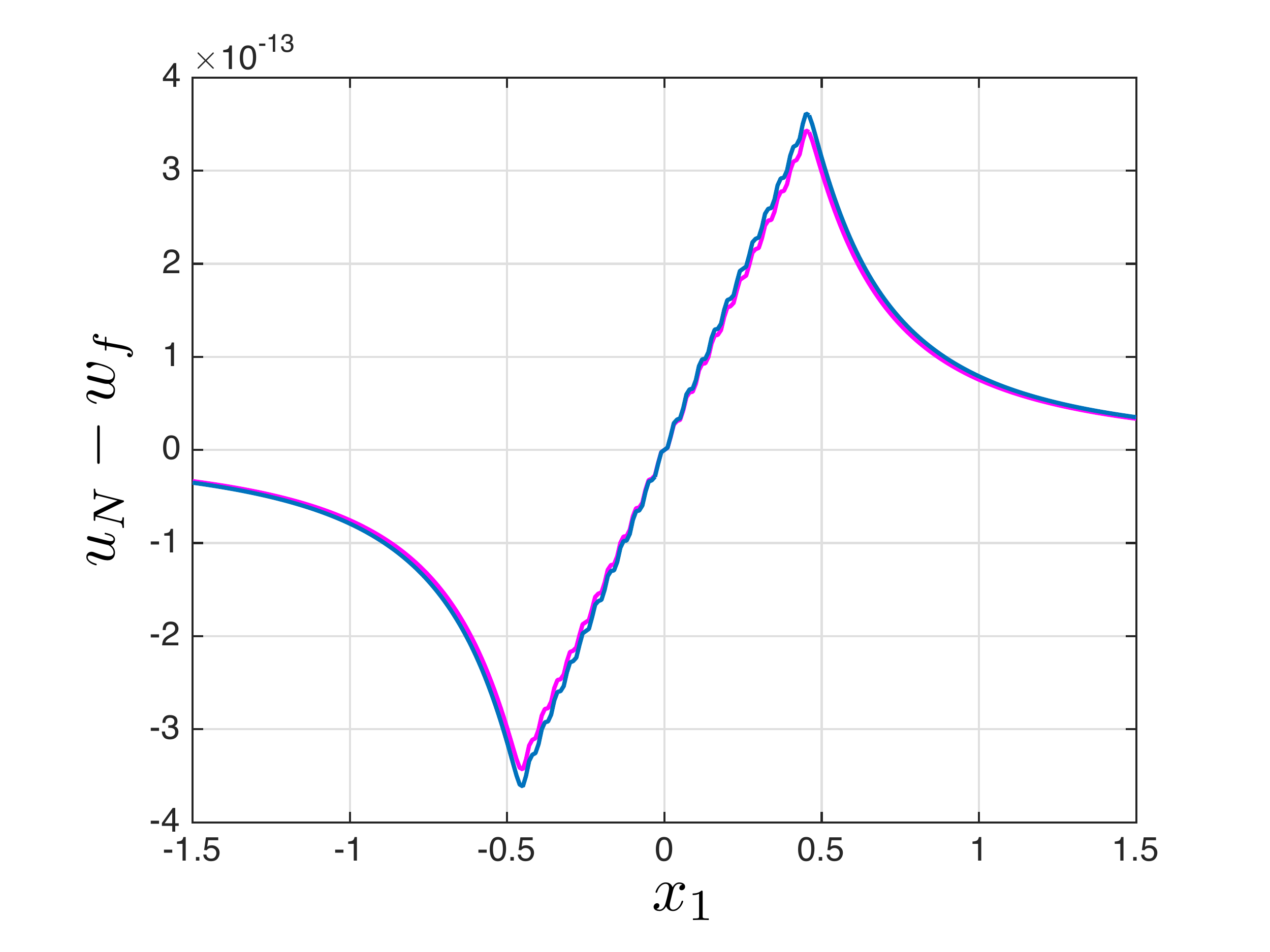}\label{fig:6b}}\\

\caption[]{Approximation to $u_N-w_f$ plotted for $-1.5\le x_1\le 1.5$, $x_2=x_3=0$. In (a), this term is shown for $N=304, 2284$ and  $5880$, ($N_1=1000, 5832$ and $13284$, respectively). In addition in (a) we supply the plot of $\hat{u}-w_f$ (see (\ref{hatsol1}), (\ref{hatsol2}) and the auxiliary homogenised problem (\ref{probuh})), which corresponds to the case ``$N\to \infty$" and  is shown by the red curve. In (b), we plot the asymptotic approximation (\ref{introeq1inf}), (\ref{alg_s_introinf}) to  $u_N-w_f$ (blue line) and  for comparison we show (\ref{introeq1inf}) where $\BC^{(j)}$, $1\le j \le N$, are computed using (\ref{Ch}) (pink line). The computations in (b) are carried out for the case  $N=5880$ ($N_1=13284$).}
\label{fig:6}%
\end{figure}

We set $x_2=x_3=0$,  $\beta=0.09$ (see (\ref{epsdnum})) and for various values of $N$ plot the asymptotic approximation for $u_N-w_f$  (see (\ref{uninf1})--(\ref{uninf5})), using (\ref{introeq1inf}) as a function of $x_1$, where $-1.5\le x_1\le 1.5$.  The line defined by $-1.5\le x_1\le 1.5$, $x_2=x_3=0$ passes  through the spherical cluster $\omega$ described in section \ref{geomspherecluster},  but does not intersect any of the inclusions. If $x_1$ increases, we see in Figure \ref{fig:6a} that as we pass through $\omega$ ($-0.5\le x_1\le 0.5$), and in particular the origin, the field undergoes a change in sign. Moreover, between $-0.5 \le x_1 \le 0.5$, one can see that the field oscillates and the number of oscillations depends on the number of inclusions in the cloud, whereas outside this region the field $u_N-w_f$ decays as is expected. 

The function $\hat{u}-w_f$ defined by (\ref{hatsol1}), (\ref{hatsol2})  (see section \ref{auxprobexp})  is also shown in Figure \ref{fig:6a}. Note this field does not oscillate inside the region $-0.5 \le x_1 \le 0.5$, and does not take into account the presence of individual inclusions. It is apparent that as $N$ increases, we see the term $u_n-w_f$ converges to the function $\hat{u}-w_f$.

On the other hand, as mentioned before, the solution to problem (\ref{probuh}) is useful in that it provides an approximation for the coefficients $\BC^{(j)}$ when $N$ is  large and can be used in the asymptotic approximation (\ref{introeq1inf}), as opposed to solving the algebraic system (\ref{alg_s_introinf}) of size $3N\times 3N$,  which can be computationally intensive.

Indeed, using (\ref{Ch}) in place of $\BC^{(j)}$ in  (\ref{introeq1inf}), for $N=5880$ ($N1=13824$) we plot the term $u_N-w_f$ in Figure \ref{fig:6b}. It is observed that the resulting plot agrees very well with results based on (\ref{introeq1inf}),  where the coefficients are computed from solving the system (\ref{alg_s_introinf}).
The procedure demonstrated here, works well in the case when periodicity is prevalent in the cluster. 
For non-periodic clusters, the solution to the auxiliary homogenised problem cannot be used to calculate $\BC^{(j)}$, $1\le j \le N$. However, in this case, the asymptotic approximation (\ref{introeq1inf}) with the coefficients are determined from (\ref{alg_s_introinf}) can handle this situation and takes into account a variety of small inclusions, whose shape and size could be different, along with the material inside each inclusion.

\section{Conclusions and discussion}\label{conclusions}
Here, we have constructed and justified a uniform asymptotic approximation for the solution to a transmission problem in a body containing many inclusions. The approximation contains a term which characterises the mutual interaction of the inclusions of within the cluster. This term makes use of the dipole fields   of individual inclusions and weights which are solutions of an algebraic system involving and integral characteristics for the inclusions. Such a term does not appear in  approximations for dilute composites. We note the approximation developed here  serves the cases when there is a dilute collection of inclusions  and when there are non-periodic arrangements of arbitrary small inclusions in a body, with $\varepsilon < \text{const } d$, for a sufficiently small constant.
 
The approximation has been shown to agree well with independent finite element computations in COMSOL and is capable of producing an accurate approximation of the solution to the transmission problem for a body with a large cluster, where finite element packages may struggle. In particular, the algebraic system governing the interaction of the inclusions within the  cluster has been linked to an auxiliary homogenised problem for an effective inclusion in a body, that relies on the cloud of inclusions being periodic. A solution of such a homogenised  problem does not take into account the oscillatory behaviour of the field in the vicinity of the defects in the cluster. However, this solution provides an effective alternative to the computation of the weights from the algebraic system, when the number of periodically placed inclusions is large. For large non-periodic arrangements of inclusions, this approach is not applicable, but the asymptotic formula constructed here remains efficient in this case.

\renewcommand{\theequation}{A.\arabic{equation}}
\setcounter{equation}{0}
\renewcommand{\thesection}{A}
\setcounter{subsection}{0}
\renewcommand{\thesubsection}{A.\arabic{subsection}}
\setcounter{subsection}{0}

\section*{Appendix: Proofs of auxiliary results}

Before presenting the proof of Lemma \ref{lemest_innerprod} (see section \ref{algmesotran}), we introduce an additional auxiliary result which we shall use.

First we  introduce the piece-wise constant functions $\Xi$ and $\Theta$ as
\begin{equation}\label{defxi}
\Xi(\Bx)=\left\{ \begin{array}{ll}
\BCP^{(j)}_\varepsilon\BC^{(j)}\;,& \quad \text{ for }\Bx \in \overline{B^{(j)}_{d/4}}\;, j=1,\dots,N\;,\\
\BO\;,& \quad \text{ otherwise}, 
\end{array}
\right.
\end{equation}
and 
\begin{equation}\label{defTheta}
\Theta(\Bx)=\left\{ \begin{array}{ll}
\BCQ^{(j)}_\varepsilon\BC^{(j)}\;,& \quad \text{ for }\Bx \in \overline{B^{(j)}_{d/4}}\;, j=1,\dots,N\;,\\
\BO\;,& \quad \text{ otherwise}.    
\end{array}
\right.
\end{equation}

We have the next result.

\begin{lem}\label{lemidentity}
The identity
\begin{equation}
\sum_{1\le i,j\le 3}\int_\Omega  \Theta_i(\BZ)\frac{\partial }{\partial Z_i} \int_\Omega  \frac{\partial }{\partial W_j} G(\BZ,\BW)\Xi_j(\BW)\, d\BW d\BZ=0\;.\label{ineqsolv3}
\end{equation}
holds.
\end{lem}
\emph{Proof. }Define 
\begin{equation}\label{feq1}
g(\BZ)=\sum_{1\le j\le 3}\int_\Omega  \frac{\partial }{\partial W_j} G(\BZ,\BW)\Xi_j(\BW)\, d\BW\;,
\end{equation}
so that the integral in the left-hand side of (\ref{ineqsolv3}) becomes
\begin{equation}
\sum_{1\le i\le 3}\int_\Omega  \Theta_i(\BZ)\frac{\partial g(\BZ) }{\partial Z_i} d\BZ\;.\label{ineqsolv31}
\end{equation}
The  function  $g(\BZ)$ of (\ref{feq1}) satisfies the problem
\[g(\BZ)=0\;, \quad \BZ\in \partial \Omega\;. \]
We apply Laplace's operator to $g$ (see (\ref{feq1})) in $\Omega$, to give
\begin{eqnarray*}
\Delta_{\BZ}g(\BZ)&=&-\sum_{1\le j \le 3}\int_\Omega \Xi_j(\BW)  \frac{\partial }{\partial W_j} (\delta(\BZ-\BW))\, d\BW\\
&=&\sum_{1\le k \le N} \sum_{1\le j\le 3}\int_{B^{(k)}_{d/4}}(\BCP_\varepsilon^{(k)}\BC^{(k)})_j  \frac{\partial }{\partial W_j} (\delta(\BW-\BZ))\, d\BW
\end{eqnarray*}
where the definitions of $G$ and $\Xi$ have be implemented in the derivation of the last result.
Next, it remains to apply integration by parts  inside $B^{(k)}_{d/4}$ to the integrals in the above right-hand side and consider $\BZ\in \Omega$. Thus, $\Delta g(\BZ)=0$ almost everywhere in $\Omega$ and using Green's representation for the function $g(\BZ)$ we deduce $g(\BZ)=0$, $\BZ\in \Omega$. Further,  consulting (\ref{ineqsolv31}) then  gives (\ref{ineqsolv3}).
Thus, the proof of Lemma \ref{lemidentity} is complete. \hfill $\Box$

\subsection*{Proof of Lemma \ref{lemest_innerprod}}
The inner product $\langle\BT\BP_\varepsilon \BCC,\BQ_\varepsilon \BCC\rangle$ appearing in (\ref{scalareq}) can be written as
\begin{equation}\label{innerprod}
\langle\BT\BP_\varepsilon \BCC,\BQ_\varepsilon \BCC\rangle=\sum_{1 \le j\le N} (\BCQ_\varepsilon^{(j)} \BC^{(j)})^T \sum_{\substack{ k \ne j \\ 1 \le k \le N}}(\nabla_{\Bz} \otimes \nabla_{\Bw}) G(\Bz, \Bw)\Big|_{\substack{\Bz=\BO^{(j)}\\ \Bw=\BO^{(k)}}} (\BCP^{(k)}_\varepsilon \BC^{(k)}) \;.
\end{equation}

The   mean value theorem for harmonic functions leads to 
\begin{eqnarray*}
(\nabla_{\Bz} \otimes \nabla_{\Bw}) G(\Bz, \Bw)\Big|_{\substack{\Bz=\BO^{(j)}\\ \Bw=\BO^{(k)}}}&=&\frac{48}{\pi d^3} \int_{B^{(k)}_{d/4}} (\nabla_{\Bz} \otimes \nabla_{\BW})G(\Bz,\BW)\Big|_{\Bz=\BO^{(j)}}\, d\BW \;.
\end{eqnarray*}
Placing this inside the inner sum of (\ref{innerprod}) gives
\begin{eqnarray*}
\langle\BT\BP_\varepsilon \BCC,\BQ_\varepsilon \BCC\rangle&=&\frac{48}{\pi d^3} \sum_{1 \le j\le N} (\BCQ_\varepsilon^{(j)} \BC^{(j)})^T \sum_{\substack{ k\ne j \\1 \le k \le N }} \int_{B^{(k)}_{d/4}} (\nabla_{\Bz} \otimes \nabla_{\BW})G(\Bz,\BW)\Big|_{\Bz=\BO^{(j)}}\, d\BW\, (\BCP^{(k)}_\varepsilon \BC^{(k)})\;.
\end{eqnarray*}
A second application of the mean value theorem  then yields:
\begin{eqnarray}
\langle\BT\BP_\varepsilon \BCC,\BQ_\varepsilon \BCC\rangle&=&\frac{48^2}{\pi^2 d^6} \sum_{1 \le j\le N} \sum_{1\le k\le N}(\BCQ_\varepsilon^{(j)} \BC^{(j)})^T  \int_{B^{(j)}_{d/4}}\int_{B^{(k)}_{d/4}} (\nabla_{\BZ} \otimes \nabla_{\BW})G(\BZ,\BW)\, d\BW d\BZ \,(\BCP^{(k)}_\varepsilon \BC^{(k)}) \nonumber\\
&&-\frac{48^2}{\pi^2 d^6} \sum_{1 \le j\le N}(\BCQ_\varepsilon^{(j)} \BC^{(j)})^T  \int_{B^{(j)}_{d/4}} \int_{B^{(j)}_{d/4}} (\nabla_{\BZ} \otimes \nabla_{\BW})G(\BZ,\BW)\, d\BW d\BZ(\BCP^{(j)}_\varepsilon \BC^{(j)})\;.\label{doublesum1}
\end{eqnarray}
Integration by parts shows that
\[ \int_{B^{(j)}_{d/4}} (\nabla_{\BZ} \otimes \nabla_{\BW})G(\BZ,\BW)\, d\BW=\int_{\partial B_{d/4}^{(j)}} (\Bn^{(j)} \otimes \nabla_{\BZ})^TG(\BZ,\BW)\, dS_{\BW}\;,\]
where $\Bn^{(j)}$ is the unit-outward normal to $B^{(j)}_{d/4}$. Here, both expressions either side of the above equation are  harmonic for $\BZ\in B_{d/4}^{(j)}$, and as a result (\ref{doublesum1}), due to the mean value theorem, becomes
\begin{eqnarray}
\langle\BT\BP_\varepsilon \BCC,\BQ_\varepsilon \BCC\rangle&=&\frac{48^2}{\pi^2 d^6} \sum_{1 \le j\le N} \sum_{1\le k\le N}(\BCQ_\varepsilon^{(j)} \BC^{(j)})^T  \int_{B^{(j)}_{d/4}}\int_{B^{(k)}_{d/4}} (\nabla_{\BZ} \otimes \nabla_{\BW})G(\BZ,\BW)\, d\BW d\BZ \,(\BCP^{(k)}_\varepsilon \BC^{(k)}) \nonumber\\
&&-\frac{48}{\pi d^3} \sum_{1 \le j\le N}(\BCQ_\varepsilon^{(j)} \BC^{(j)})^T   \int_{\partial B^{(j)}_{d/4}} (\nabla_{\Bz} \otimes \nabla_{\BW})G(\Bz,\BW)\Big|_{\Bz=\BO^{(j)}}\, d\BW (\BCP^{(j)}_\varepsilon \BC^{(j)})\;.\label{doublesum}
\end{eqnarray}
 The fact $G(\Bx, \By)=O(|\Bx-\By|^{-1})$ allows one to derive the inequality
\[\Big|\int_{\partial B^{(j)}_{d/4}} (\Bn^{(j)} \otimes \nabla_{\BW})G(\Bz,\BW)\Big|_{\Bz=\BO^{(j)}}\, d\BW\Big|\le \text{Const}\;.\]
Thus application of the Cauchy inequality and the preceding inequality shows that
\begin{eqnarray}
&& \sum_{1 \le j\le N}(\BCQ_\varepsilon^{(j)} \BC^{(j)})^T  \int_{\partial B^{(j)}_{d/4}} (\nabla_{\Bz} \otimes \nabla_{\BW})G(\Bz,\BW)\Big|_{\Bz=\BO^{(j)}}\, d\BW (\BCP^{(j)}_\varepsilon \BC^{(j)})\nonumber\\
&\le &\text{Const }\sum_{1\le j \le N}|\BCQ^{(j)}_\varepsilon \BC^{(j)}||\BCP^{(j)}_\varepsilon \BC^{(j)}| \nonumber\\
&\le &\text{Const }\sum_{1\le j \le N}|\BCQ^{(j)}_\varepsilon \BC^{(j)}|^2\;,\label{ineqsolv1}
\end{eqnarray}
where in the last step the definition of $\BCQ_\varepsilon^{(j)}$ in (\ref{eqDeps}) has also been used.

The definitions of $\Theta$ and $\Xi$ in (\ref{defTheta}) and (\ref{defxi}), respectively, allow the double sum in (\ref{doublesum}) to take the equivalent form
\begin{eqnarray}
&&\sum_{1 \le j\le N} \sum_{1\le k\le N}(\BCQ_\varepsilon^{(j)} \BC^{(j)})^T  \int_{B^{(j)}_{d/4}}\int_{B^{(k)}_{d/4}} (\nabla_{\BZ} \otimes \nabla_{\BW})G(\BZ,\BW)\, d\BW d\BZ \,(\BCP^{(k)}_\varepsilon \BC^{(k)})
\nonumber \\
&&=\sum_{i,j=1}^3\int_\Omega  \Theta_i(\BZ)\frac{\partial }{\partial Z_i} \int_\Omega  \frac{\partial }{\partial W_j} G(\BZ,\BW)\Xi_j(\BW)\, d\BW d\BZ\;.\label{ineqsolv2}
\end{eqnarray}

Combining (\ref{innerprod})--(\ref{ineqsolv1}), (\ref{ineqsolv2}) and (\ref{ineqsolv3}) shows that 
\[|\langle\BT\BP_\varepsilon \BCC,\BQ_\varepsilon \BCC\rangle| \le \text{Const } d^{-3}\langle\BQ_\varepsilon \BCC, \BQ_\varepsilon \BCC\rangle \;,\]
completing the proof of Lemma \ref{lemest_innerprod}. \hfill $\Box$
 
\end{document}